\documentclass[12pt]{amsart}
\usepackage{amsmath,amsfonts,amssymb,amscd}
\usepackage{latexsym}

\usepackage[
hyperindex=true,pagebackref=true,bookmarks=true,
colorlinks=true,linkcolor=blue,citecolor=red]
{hyperref}
\def\~{{\rm --}} 

\title [Whittaker limits of difference spherical functions]
{Whittaker limits of difference spherical functions} 
\author[Ivan Cherednik]{Ivan Cherednik $^\dag$}

\thanks{$^\dag$  \today\ \ \ Partially supported by NSF grant
DMS--0800642}

\address[I. Cherednik]{Department of Mathematics, UNC
Chapel Hill, North Carolina 27599, USA\\
chered@email.unc.edu}

 \def\bysame{{\bf --- }}
 \def\~{{\bf --}}
\newcommand{\comment}[1]{}
\renewcommand{\tilde}{\widetilde}
\renewcommand{\hat}{\widehat}
\newcommand{\dagx}{\hbox{\tiny\mathversion{bold}$\dag$}}
\newcommand{\ddagx}{\hbox{\tiny\mathversion{bold}$\ddag$}}

\renewcommand{\tilde}{\widetilde}
\renewcommand{\hat}{\widehat}

\newcommand{\Z}{{\mathbb Z}}
\newcommand{\Q}{{\mathbb Q}}
\newcommand{\N}{{\mathbb N}}
\newcommand{\C}{{\mathbb C}}
\newcommand{\R}{{\mathbb R}}

\def\HH{\mbox{${\mathcal H}$\kern-5.2pt${\mathcal H}$}}

\newtheorem{theorem}{Theorem}[section]

\newtheorem{proposition}[theorem]{Proposition}
\newtheorem{definition}[theorem]{Definition}
\newtheorem{lemma}[theorem]{Lemma}
\newtheorem{corollary}[theorem]{Corollary}

\newtheorem{theorem }{Theorem}[section]
\newtheorem{maintheorem }[theorem]{Main Theorem}
\newtheorem{proposition }[theorem]{Proposition}
\newtheorem{definition }[theorem]{Definition}
\newtheorem{lemma }[theorem]{Lemma}
\newtheorem{corollary }[theorem]{Corollary}
\newtheorem{notation }[theorem]{Notation}
\newtheorem{remark }[theorem]{Remark}
\newtheorem{example }[theorem]{Example}

\newtheorem{ maintheorem }[theorem]{Main Theorem}
\newtheorem{ theorem}{Theorem}[section]
\newtheorem{ proposition}[theorem]{Proposition}
\newtheorem{ definition}[theorem]{Definition}
\newtheorem{ lemma}[theorem]{Lemma}
\newtheorem{ corollary}[theorem]{Corollary}
\newtheorem{ notation}[theorem]{Notation}
\newtheorem{ remark}[theorem]{Remark}
\newtheorem{ example}[theorem]{Example}

\def\for{\  \hbox{ for } \ }
\def\iif{ \ \hbox{ if } \ }
\def\when{ \ \hbox{ when } \ }
\def\where{\  \hbox{ where } \ }
\def\and{\  \hbox{ and } \ }
\def\and{\  \hbox{ and } \ }

\def\equal{\stackrel{\,\mathbf{def}}{= \kern-3pt =}}

\def\la{\lambda}
\def\La{\Lambda}
\def\om{\omega}

\def\Th{\Theta}
\def\th{\theta}
\def\al{\alpha}

\def\ga{\gamma}

\def\de{\delta}

\def\si{\sigma}

\def\Ga{\Gamma}
\def\ze{\zeta}

\def\vth{{\vartheta}}

\def\tal{\tilde{\alpha}}

\def\tga{\tilde{\gamma}}
\def\tGa{\tilde{\Gamma}}

\def\tw{\widetilde w}
\def\tW{\widetilde W}

\def\tz{\tilde z}
\def\tb{\tilde b}

\def\tR{\tilde R}

\def\hw{\widehat{w}}
\def\hW{\widehat{W}}
\def\hu{\hat{u}}

\def\hv{\hat{v}}

\def\P{\mathbf{P}}

\def\0{\mathbf{0}}

\def\çF{\mathcal{F}}

\def\r{\mathcal{R}}

\def\p{\mathcal{P}}

\def\h{\mathcal{H}}

\def\e{\mathcal{E}}
\def\v{\mathcal{V}}

\def\lan{\langle}
\def\llb{(\!(}
\def\ran{\rangle}
\def\rrb{)\!)}
\def\lng{\hbox{\rm{\tiny lng}}}
\def\sht{\hbox{\rm{\tiny sht}}}

\newcommand{\sq}{\phantom{1}\hfill$\qed$}

\def\HH{\mathfrak{H}}

\def\HH{\hbox{${\mathcal H}$\kern-5.2pt${\mathcal H}$}}

\font\smm=msbm10 at 12pt 
\def\symbol#1{\hbox{\smm #1}}
\def\lsmash{{\symbol n}}

\def\#{\sharp}

\begin{document}
\maketitle
\renewcommand{\baselinestretch}{1.2}. 
{\textmd
\tableofcontents
} 
\renewcommand{\baselinestretch}{1.0}.
\vfill\eject

\renewcommand{\natural}{\wr}

\setcounter{section}{-1}
\setcounter{equation}{0}
\section{Introduction}

\comment{
We introduce the q-Whittaker functions as the limit at t=0 
of the q,t-spherical function extending the symmetric 
Macdonald polynomials to arbitrary eigenvalues. The 
construction heavily depends on the technique of the 
q-Gaussians developed by the author (and Stokman in the
non-reduced case). In this approach, the q-Whittaker function 
is given by a series convergent everywhere. One of the 
applications is a q-version of the Shintani-Casselman-Shalika 
formula, which is directly connected with q-Mehta-Macdonald 
identities in terms of the Jackson integrals. This formula 
generalizes that of type A due to Gerasimov et al. to 
arbitrary reduced root systems. At the end of the paper,
we generalize the Harish-Chandra asymptotic formula 
for the spherical functions to the q,t-case, including the
Whittaker limit.
}
The main aim of this paper is to introduce {\em global}
\, $q$\~Whittaker functions as the limit $t\to 0$ of
the (renormalized) generalized symmetric 
spherical functions constructed in \cite{C5}
for arbitrary reduced root systems (see \cite{Sto} in the
$C^\vee C$\~case). This work is inspired by \cite{GLO1} and
\cite{GLO2}, though our approach is different. 
For instance, we obtain a $q$\~version of the classical 
Shintani-Casselman-Shalika formula \cite{Shi,CS} 
via the\, $q$\~Mehta -Macdonald 
integral in the Jackson setting. The 
Shintani-type formulas (in the case of $GL_n$) play
an important role in \cite{GLO1,GLO2}, but 
the\, $q$\~Gauss integrals are not considered there
as well as globally-defined\, $q$\~Whittaker functions.
We use these formulas to obtain a $q,t$\~generalization
of the Harish-Chandra  
asymptotic formula for the classical spherical function.
\smallskip

\subsection{Results and applications}
The key observation is that
the definition of the symmetric $q,t$\~spherical functions
from \cite{C5} is compatible with taking the Whittaker
limit and results in globally-defined \,
$q$\~Whittaker functions. The definition from \cite{C5} 
is based on the\, $q$\~Mehta- Macdonald integrals calculated 
there for the constant term functional, i.e., in the setting of
Laurent series. In this paper, we mainly treat the spherical
functions as {\em global} ones, analytic or meromorphic.
 
The $q$\~Whittaker functions are solutions of the $q$\~Toda
eigenvalue problem and are expected to have 
important applications in mathematics 
and physics, including the Langlands program. Concerning
the latter and relations to the affine flag varieties, see, 
for instance, \cite{GiL,BF,Ion2}. 
The\, $q$\~Shintani-Casselman-Shalika formula gives a (relatively
simple) example of the Langlands correspondence.
The affine Toda lattice provides another link;
it is (presumably) dual to the $q$\~Toda lattice in the
sense of \cite{KL2} (via the monodromy map).

The coefficients of the expansion
of our Whittaker function are essentially
polynomials in terms of $q$
with {\em non-negative integral} coefficients. It can be verified
using the intertwining operators 
or via the relation to the Demazure characters (we do not
discuss it in this paper). 
This fact is of obvious importance for the ``categorization" of the 
$q$\~Whittaker function and its geometric applications. 
The ``$q$\~integrality" has no known counterpart
in the general $q,t$\~theory (with a reservation concerning
the {\em stable} $GL$\~case); one of the parameters, 
$q$ or $t$, has to be eliminated or expressed in terms of 
the remaining one. However, the $q,t$\~spherical functions 
are more convenient to deal with in many other aspects,
including the analytic theory. 

The cornerstone of their theory is the
{\em duality} based on the DAHA-Fourier transform;
see \cite{C101} and\cite{C103}. It is a special feature
of the general $q,t$\~theory, missing for
the $q$\~Whittaker functions and 
$q$\~Hermite polynomials ($t\to 0$),
in the Harish-Chandra theory ($q\to 1$) and in 
the $p$\~adic limit ($q\to \infty$).
A specific problem with the Whittaker limiting procedure 
among other degenerations is that it destroys the 
$W$\~invariance; it gives another reason for
treating Whittaker functions
as limits of the spherical functions
rather than for creating their intrinsic theory.
On the other hand, the Whittaker functions satisfy quite
a few identities that are not present in the 
$q,t$\~theory. These formulas, the $q$\~integrality
of the coefficients and various applications 
obviously make the $q$\~Whittaker functions an important
independent direction, which requires developing specific 
methods.

At the end of the paper, we outline the approach to
the global spherical and Whittaker functions via the
harmonic analysis. Our formulas for these functions
are actually equivalent to certain fundamental properties 
of the corresponding integral transforms in the space of 
Laurent polynomials multiplied by the Gaussian. The
latter space is the simplest and the most natural choice here,
but the same functions can serve other algebraic and analytic
situations. The harmonic analysis direction
seems very promising. For
instance, the existence of the $q$\~Whittaker
limit of the global spherical function appears a {\em boundary}
case of the general theory of {\em growth estimates} 
for the $q,t$\~spherical function in terms of $x,\la$(the
spectral parameter) and $t=q^k$.
\smallskip

\subsection{Growth estimates}
Provided that $\Re (x), \Re(\la)$ are inside the positive 
Weyl chamber $\mathfrak{C}_+$ 
(the walls must be avoided), 
the global spherical function for $0<q<1$
approaches asymptotically in the limit of 
large $\Re (x)$ 
$$
|W\!|\ \hbox{C\!T} \ \Theta(\rho_k)\ 
\frac{\Theta(x+\la-\rho_k)}
{\Theta(x)\Theta(\la)}\, \prod_{\al\in R_+}
\frac{\Gamma_q(\la_\al^\vee)}
{\Gamma_q(\la_\al^\vee+k_\al)}
$$
for the theta-series $\Th$ and $q$\~Gamma function
associated with a given root system $R$;\, $t=q^k$,\,  
$\rho_k=k\rho$\, in the simply-laced case,\, C\!T being
the constant term of the celebrated (symmetric) Macdonald 
function. No inequalities
for $k$ are necessary but one must avoid the values 
where the polynomial representation of DAHA  
becomes non-semisimple. 

Up to a periodic function, the $x$\~dependence of this
function is $q^{(x,\,\rho_k-\la)}$, so this theorem
is an {\em exact} $q,t$\~analog
of the Harish-Chandra formula \cite{HC}
describing the asymptotic behavior of the classical
spherical function in terms of the $c$\~function. 
In the Whittaker limit, 
$\rho_k$ is omitted and  $\Re (x)$ must be taken from 
$-\mathfrak{C}_+$  (the Whittaker function
it is not $W$\~invariant with respect to $x$).

It is one of the major results of this
paper, which seems a beginning of fruitful 
analytic $q$\~theory. 
\smallskip
   
\subsection{Our approach}
It is different from that of \cite{GLO1,GLO2}
(and we deal with arbitrary reduced root systems).  
The technique of the Gaussians is the key to introduce
the {\em global}\, $q$\~Whittaker function and prove 
the Shintani-type formulas. The $q$\~Whittaker function 
is mainly treated in \cite{GLO1,GLO2} as a {\em discrete} 
function on the weight lattice for $GL_n$ satisfying 
the $q$\~Toda system of {\em difference} equations. 

The space of all solutions is, generally,
$|W\!|$\~dimensional over the field of periodic functions, 
playing the role of constants in the difference theory; 
upon the restriction to the weight lattice 
it is $|W\!|$\~dimensional over $\C$. Choosing the 
``right" Whittaker function 
in this space requires certain growth conditions; using the 
$W$\~symmetric dependence on the spectral parameters gives another
approach. There is no intrinsic definition of the $q$\~Whittaker 
function so far, but our formula and the growth conditions we
establish clarify what can be expected. First, only
positive powers appear in its Laurent series expansion
(after dropping the Gaussians). Second, our $x$\~asymptotic
formula for the $q$\~Whittaker function inside the negative
Weyl chamber is sufficient to fix it uniquely.

We note that in the {\em differential} setting, 
the spherical and Whittaker functions can be uniquely
determined from the eigenvalue problem (subject to
the $W$\~invariance for the spherical function and certain 
growth conditions in the Whittaker case). It simplifies 
the starting definitions.  However the difference theory 
is more universal and, remarkably, has 
important algebraic and analytic advantages. The self-duality 
of the DAHA-Fourier transform and the technique of the Gaussians
are the key; these are special features of the $q,t$\~setting 
and are mainly absent in the trigonometric-differential 
and $p$\~adic cases. 
In this respect, the $q,t$\~theory is somewhat similar 
to the rational-differential theory of (multi-variable) 
Bessel functions.  
\smallskip

\subsection{Difference spherical functions}
The {\em global} nonsymmetric and symmetric 
{\em $q,t$\~spherical functions} were defined in 
\cite{C5} and then in \cite{Sto} (the 
$C^\vee C$\~case) as the reproducing kernels of the Fourier
transform of the standard polynomial representation twisted by 
the Gaussian. In this approach, the spherical function is
determined uniquely (the Macdonald eigenvalue problem fixes 
it only up to periodic factors). Using the Gaussians, among other
things, provides the global convergence.
These construction appeared compatible with the Whittaker limit.

The Gaussians play the key role in  
our approach to the Shintani-Casselman-Shalika 
formula. In the $q,t$\~setting, it becomes
the Mehta -Macdonald formula in the Jackson case 
from \cite{C5}, where a special vector,
$-\rho_k$, is taken as the origin of the Jackson summation. 

Developing this direction, we conclude the paper with the 
Jackson-Gauss integrals for the {\em global} spherical and 
Whittaker functions; such formulas were given only for Macdonald 
polynomials in \cite{C5}. These formulas seem an important
step toward systematic difference harmonic analysis, although
the case of the real integration is still beyond the
existing theory. Now, with the $q,t$\~Harish-Chandra asymptotic 
formulas from this paper, it seems that there are no 
obstacles for developing the real integration theory generalizing
the classical ``non-compact" case.
\smallskip

Conceptually, as it was observed in \cite{GLO2}, 
the $q$\~variant of the Shintani-Casselman-Shalika formula 
is nothing but the duality formula for the Macdonald polynomials 
from \cite{C3} considered upon the limit $t\to 0$. 
However, establishing exact relations is, generally, 
a subtle problem. The Shintani-type formulas play the major
role in the paper, including the growth estimates.

This interpretation gives evidence that the DAHA-Fourier transform 
is connected with the (local quantum) geometric Langlands 
correspondence. Generally, the DAHA\~localization functor,
which includes the modular transformation $q\mapsto q'$, 
is expected to play its role in the quantum geometric Langlands 
correspondence; the DAHA-Fourier transform is likely to be one 
of its ingredients. 
\smallskip

We note that DAHA leads to a theory
that is a priori more general than the one needed
for the (local) quantum Langlands correspondence because 
it contains an extra
parameter $t$. However, there is growing evidence that the
general $q,t$\~DAHA appear in the Langlands program.
It makes important the exact relations
between the $q,t$\~spherical functions and $q$\~Whittaker ones
(which are already a part of the Langlands program). 
We expect this paper to trigger interesting new developments. 
\smallskip

It is worth mentioning that the approach to spherical
functions via the Fourier transform 
depends on the choice of the corresponding  
representation of the double affine Hecke algebra. 
Technically, the choice of this space influences only
the normalization; spherical function are defined up
to periodic factors. However, the analytic properties
of the $q,t$\~spherical function, exact factors in  
the Shintani-type formulas and other similar
features reflect the properties of the considered
representation (equivalently, the choice of
the normalization). 

For instance, if the 
Gaussian is interpreted as a theta-function, then
the corresponding spherical function is {\em meromorphic} 
but not analytic. Treating the Gaussian as $q^{x^2/2}$ (not as 
a Laurent series), i.e., using a somewhat different analytic
setting, leads to the $q,t$\~spherical functions
analytic everywhere, but not single-valued in terms of $q^x$. 
If the Gaussians are omitted
in this definition, i.e., the DAHA-Fourier transform
acts from the polynomial representation to the space
of delta-functions, then the corresponding spherical
function will become a generalized function. A general
problem is in finding a representation that ensures the
best analytic properties of the reproducing kernel;
if $|q|\,\hbox{\tiny ${}^\gtrless$}\, 1$, then 
the polynomial representation times the Gaussian is the one.
\smallskip

\subsection{The setting of the paper}
Only the symmetric theory will be considered in this
work; the (truly) nonsymmetric $q$\~Whittaker function can be 
defined as certain limits of the {\em nonsymmetric} global
spherical function, but the construction becomes more involved and 
will be a subject of the next work(s).
Nevertheless, we begin the paper
with the account of the nonsymmetric Macdonald polynomials 
including their (straight) degeneration 
as $t\to 0$, which is  closely related
to the Demazure characters of irreducible affine Lie
algebras; see \cite{San,Ion1}.  
We mainly need the formulas in terms of the intertwining 
operators to justify some of our claims and estimates; 
the intertwiners can be naturally defined
only in the nonsymmetric theory. 

We mention that the Macdonald symmetric
polynomials considered under the limit $t\to 0$ 
generalize the classical $q$\~Hermite polynomials,
so the main result of the paper is in establishing
the formula for the $q$\~Whittaker function in terms of
{\em multi-variable\, $q$\~Hermite polynomials}.

In the theory of nonsymmetric Whittaker functions 
(it is beyond this paper and not completed so far), the 
{\em nonsymmetric} $q$\~Whittaker 
function become a {\em generating function for 
all Demazure characters}, not only the ones for
anti-dominant weights. To be exact, $q$\~Hermite
polynomials appear here instead of the Demazure characters
(there is a direct link). However, this interpretation
requires new technique of $W$\~spinors, and the analytic
aspects are not clear at the moment. The appearance of
all Demazure characters can clarify the role of Whittaker 
functions in the Kac-Moody theory and may have connections to 
\cite{GiL} (quantum $K$\~theory of affine flag varieties)
and to questions and conjectures from \cite{BF} concerning
the $IC$\~theory of affine flag varieties.  

We note that there are two possible setups in the DAHA theory 
for the non-simply-laced root systems, which
correspond to two possible choices of the affine extension.
In this paper, we introduce the affine root 
system using $\al_0=[-\vth,1]$
in terms of the maximal {\em short} root $\vth$ (the
so-called twisted case).
The conjugation by the Gaussian and the Fourier transform 
preserve the double affine Hecke algebra in this setup. 
By the way, it is exactly the case when a relation
to the Demazure characters can be established according to
\cite{Ion1}, Theorem 1.

The case of the ``standard" non-twisted affine root system with 
$\al_0=[-\th,1]$ for the maximal {\em long} root $\th$ is
analogous, although the Fourier transform acts from the 
double affine Hecke algebra to its dual in the $B,C$\~cases.
This setting is expected to be related to the geometric Langlands 
correspondence (cf. \cite{C103}).

Technically, the switch to the standard
DAHA can be achieved by changing the action of $T_0$ in the 
polynomial representation. This change influences the relations 
of $T_0$ with the $X$\~operators (indexed by the weights).
The $Y$\~operators become labeled by the coweights 
for such choice of $T_0$; they are labeled by the 
weights in this paper. However, this transformation is far from 
being direct at level of difference Mehta-Macdonald formulas
we need for the theory of difference spherical and Whittaker
functions. 


\medskip
{\bf Acknowledgements.}
The author is thankful to D.~Kazhdan for alerting
me to the works of Gerasimov et. al and for our
various conversations on the Whittaker functions
and the Langlands correspondence. I indebted
to D.~Gaitsgory for the discussion of 
the quantum geometric Langlands duality. Special thanks go to 
A.~Gerasimov for his explanations of the results of
\cite{GLO1,GLO2}, which influenced this paper a great 
deal. I am very grateful to E.~Opdam and J.~Stokman for
reading the paper and suggesting various improvements.
 
\medskip
\setcounter{equation}{0}
\section{Double Hecke algebra}
Let $R=\{\al\}   \subset \R^n$ be a root system of type
$A,B,...,F,G$
with respect to a euclidean form $(z,z')$ on $\R^n
\ni z,z'$,
$W$ the {\em Weyl group} 
generated by the reflections $s_\al$,
$R_{+}$ the set of positive  roots ($R_-=-R_+$)
corresponding to fixed simple 
roots $\al_1,...,\al_n,$
$\Ga$ the Dynkin diagram
with $\{\al_i, 1 \le i \le n\}$ as the vertices.
Accordingly,
$$R^\vee=\{\al^\vee =2\al/(\al,\al)\}.$$

The root lattice and the weight lattice are:
\begin{align}
& Q=\oplus^n_{i=1}\Z \al_i \subset P=\oplus^n_{i=1}\Z \om_i,
\notag
\end{align}
where $\{\om_i\}$ are fundamental weights:
$ (\om_i,\al_j^\vee)=\de_{ij}$ for the
simple coroots $\al_i^\vee.$
Replacing $\Z$ by $\Z_{\pm}=\{m\in\Z, \pm m\ge 0\}$ we obtain
$Q_\pm, P_\pm.$
Here and further see  \cite{Bo}.

The form will be normalized
by the condition  $(\al,\al)=2$ for the
{\em short} roots in this paper. 
Thus,

\centerline{
$\nu_\al\equal (\al,\al)/2$ can be either $1,$ or $\{1,2\},$ or
$\{1,3\}.$ }
\noindent

This normalization leads to the inclusions
$Q\subset Q^\vee,  P\subset P^\vee,$ where $P^\vee$ is
defined to be generated by
the fundamental coweights $\{\om_i^\vee\}$ dual to $\{\al_i\}$. 
\smallskip

We set 
$\nu_i\ =\ \nu_{\al_i}, \
\nu_R\ = \{\nu_{\al}, \al\in R\}$ and
\begin{align}\label{partialrho}
&\rho_\nu\equal (1/2)\sum_{\nu_{\al}=\nu} \al \ =
\ \sum_{\nu_i=\nu}  \om_i, \hbox{\ where\ } \al\in R_+,
\ \nu\in\nu_R.
\end{align}
Note that $(\rho_\nu,\al_i^\vee)=1$ for $\nu_i=\nu.$

\subsection{Affine Weyl group}
The vectors $\ \tal=[\al,\nu_\al j] \in
\R^n\times \R \subset \R^{n+1}$
for $\al \in R, j \in \Z $ form the
{\em affine root system}
$\tR \supset R$ ($z\in \R^n$ are identified with $ [z,0]$).
We add $\al_0 \equal [-\vth,1]$ to the simple
roots for the {\em maximal short root} $\vth\in R_+$.
It is also the {\em maximal positive
coroot} because of the choice of normalization. 

The corresponding set
$\tR_+$ of positive roots equals 
$R_+\cup \{[\al,\nu_\al j],\ \al\in R, \ j > 0\}$.
Indeed, any positive affine root $[\al,\nu_\al j]$
is a linear combinations
with non-negative integral coefficients of
$\{\al_i,\,0\le i\le n\}$.

We complete the Dynkin diagram $\Ga$ of $R$
by $\al_0$ (by $-\vth$, to be more
exact); it is called {\em affine Dynkin diagram}
$\tGa$. One can obtain it from the
completed Dynkin diagram from \cite{Bo} for 
the {\em dual system}
$R^\vee$ by reversing all arrows. 

The set of
the indices of the images of $\al_0$ by all
the automorphisms of $\tGa$ will be denoted by $O$
($O=\{0\} \for E_8,F_4,G_2$). Let $O'=\{r\in O, r\neq 0\}$.
The elements $\om_r$ for $r\in O'$ are the so-called minuscule
weights: $(\om_r,\al^\vee)\le 1$ for
$\al \in R_+$.

Given $\tal=[\al,\nu_\al j]\in \tR,  \ b \in P$, let
\begin{align}
&s_{\tal}(\tz)\ =\  \tz-(z,\al^\vee)\tal,\
\ b'(\tz)\ =\ [z,\ze-(z,b)]
\label{ondon}
\end{align}
for $\tz=[z,\ze] \in \R^{n+1}$.

The 
{\em affine Weyl group} $\tW$ is generated by all $s_{\tal}$
(we write $\tW = \lan s_{\tal}, \tal\in \tR_+\ran)$. One can take
the simple reflections $s_i=s_{\al_i}\ (0 \le i \le n)$
as its
generators and introduce the corresponding notion of the
length. This group is
the semidirect product $W\lsmash Q'$ of
its subgroups $W=$ $\lan s_\al,
\al \in R_+\ran$ and $Q'=\{a', a\in Q\}$, where
\begin{align}
& \al'=\ s_{\al}s_{[\al,\,\nu_{\al}]}=\
s_{[-\al,\,\nu_\al]}s_{\al}\for
\al\in R.
\label{ondtwo}
\end{align}

The {\em extended Weyl group} $ \hW$ generated by $W\and P'$
(instead of $Q'$) is isomorphic to $W\lsmash P'$:
\begin{align}
&(wb')([z,\ze])\ =\ [w(z),\ze-(z,b)] \for w\in W, b\in B.
\label{ondthr}
\end{align}
From now on,  $b$ and $b',$ $P$ and $P'$ will be identified.

Given $b\in P_+$, let $w^b_0$ be the longest element
in the subgroup $W_0^{b}\subset W$ of the elements
preserving $b$. This subgroup is generated by simple
reflections. We set
\begin{align}
&u_{b} = w_0w^b_0  \in  W,\ \pi_{b} =
b( u_{b})^{-1}
\ \in \ \hW, \  u_i= u_{\om_i},\pi_i=\pi_{\om_i},
\label{xwo}
\end{align}
where $w_0$ is the longest element in $W,$
$1\le i\le n.$

The elements $\pi_r\equal\pi_{\om_r}, r \in O'$ and
$\pi_0=\hbox{id}$ leave $\tGa$ invariant
and form a group denoted by $\Pi$,
 which is isomorphic to $P/Q$ by the natural
projection $\{\om_r \mapsto \pi_r\}$. As to $\{ u_r\}$,
they preserve the set $\{-\vth,\al_i, i>0\}$.
The relations $\pi_r(\al_0)= \al_r= ( u_r)^{-1}(-\vth)$
distinguish the
indices $r \in O'$. Moreover,
\begin{align}
& \hW  = \Pi \lsmash \tW, \where
  \pi_rs_i\pi_r^{-1}  =  s_j \iif \pi_r(\al_i)=\al_j,\
 0\le j\le n.
\end{align}

We will need 
the following 
{\em affine action} of $\hW$ on
$z \in \R^n$:
\begin{align}
& (wb)\llb z \rrb \ =\ w(b+z),\ w\in W, b\in P,\notag\\
& s_{\tal}\llb z\rrb\ =\ z - ((z,\al^\vee)+j)\al,
\ \tal=[\al,\nu_\al j]\in \tR.
\label{afaction}
\end{align}
For instance, $(b w)\llb 0\rrb=b$ for any $w\in W.$
The relation to the above action is given in terms of
the {\em affine pairing}
$([z,l], z'+d)\equal (z,z')+l:$
\begin{align}
& (\hw([z,l]),\hw\llb z' \rrb+d) \ =\
([z,l], z'+d) \for \hw\in \hW,
\label{dform}
\end{align}
where we treat $d$ formally.
\smallskip

\subsection{The length 
\texorpdfstring{{\mathversion{bold}on $\hW$}}{}}
Setting
$\hw = \pi_r\tw \in \hW,\ \pi_r\in \Pi, \tw\in \tW,$
the length $l(\hw)$
is by definition the length of the reduced decomposition
$\tw= $ $s_{i_l}...s_{i_2} s_{i_1} $
in terms of the simple reflections
$s_i, 0\le i\le n.$ The number of  $s_{i}$
in this decomposition
such that $\nu_i=\nu$ is denoted by   $l_\nu(\hw).$

The {\em length} can be
also defined as the
cardinality $|\la(\hw)|$
of the {\em $\la$\~set} of $\hw$\,:
\begin{align}\label{lasetdef}
&\la(\hw)\equal\tR_+\cap \hw^{-1}(\tR_-)=\{\tal\in \tR_+,\
\hw(\tal)\in \tR_-\},\
\hw\in \hW.
\end{align}
Alternatively,
\begin{align}
&\la(\hw)=\cup_\nu\la_\nu(\hw),\
\la_\nu(\hw)\
\equal\ \{\tal\in \la(\hw),\nu({\tal})=\nu \}.
\label{xlambda}
\end{align}
\smallskip

The coincidence with the previous definition 
is based on the equivalence of the {\em length equality}
\begin{align}\label{ltutwa}
&(a)\ \ l_\nu(\hw\hu)=
l_\nu(\hw)+l_\nu(\hu)
\for \hw,\hu\in\hW
\end{align}
and the {\em cocycle relation}
\begin{align}
&  (b)\ \ \la_\nu(\hw\hu) = \la_\nu(\hu) \cup
\hu^{-1}(\la_\nu(\hw)),
\label{ltutw}
\end{align}
which, in its turn, is equivalent to
the {\em positivity condition}
\begin{align}\label{ltutwc}
& (c)\ \  \hu^{-1}(\la_\nu(\hw))
\subset \tR_+
\end{align}
and is also equivalent to the {\em embedding condition}
\begin{align}\label{ltutwd}
& (d)\ \  \la_\nu(\hu)\subset \la_\nu(\hw).
\end{align}

See, e.g., \cite{C4,C101} and also \cite{Bo,Hu}. 
Applying (\ref{ltutw}) to the reduced decomposition
$\hw=\pi_rs_{i_l}\cdots s_{i_2}s_{i_1},$
\begin{align}
\la(\hw) = \{\ &\tal^l=\tw^{-1}s_{i_l}(\al_{i_l}),\
\ldots,\ \tal^3=s_{i_1}s_{i_2}(\al_{i_3}),\notag\\
&\tal^2=s_{i_1}(\al_{i_2}),\ \tal^1=\al_{i_1}\  \}.
\label{tal}
\end{align}

\medskip
\subsection{Reduction modulo 
\texorpdfstring{\mathversion{bold}$W$}{\em W}}
It generalizes the construction of the elements
$\pi_{b}$ for $b\in P_+;$ see 
\cite{C4} or \cite{C101}.

\begin{proposition} \label{PIOM}
 Given $ b\in P$, there exists a unique decomposition
$b= \pi_b  u_b,$
$ u_b \in W$ satisfying one of the following equivalent conditions:

{(i) \  } $l(\pi_b)+l( u_b)\ =\ l(b)$ and
$l( u_b)$ is the greatest possible,

{(ii)\  }
$ \la(\pi_b)\cap R\ =\ \emptyset$.

The latter condition implies that
$l(\pi_b)+l(w)\ =\ l(\pi_b w)$
for any $w\in W.$ Besides, the relation $ u_b(b)
\equal b_-\in P_-=-P_+$
holds, which, in its turn,
determines $ u_b$ uniquely if one of the following equivalent
conditions is imposed:

{(iii) }
$l( u_b)$ is the smallest possible,

{(iv)\ }
if\, $\al\in \la( u_b)$  then $(\al,b)\neq 0$.
\end{proposition}
\qed

\smallskip

Condition (ii) readily
gives a complete description of the set
$\pi_P=\{\pi_b, b\in P\}$, namely, 
only
$\,[\,\al<0,\,\nu_\al j>0\,]\,$
can appear in $\la(\pi_b)$.
\smallskip

Explicitly,
\begin{align}
\la(b) = \{ \tal>0,\  &( b, \al^\vee )>j\ge 0 \iif \al\in R_+,
\label{xlambi}\\
&( b, \al^\vee )\ge j> 0 \iif \al\in R_-\},
\notag \\
\la(\pi_b) = \{ \tal>0,\ \al\in R_-,\
&( b_-, \al^\vee )>j> 0
\iif  u_b^{-1}(\al)\in R_+,
\label{lambpi} \\
&( b_-, \al^\vee )\ge j > 0 \iif
u_b^{-1}(\al)\in R_- \}, \notag 
\end{align}
For instance,
$l(b)=l(b_-)=-2(\rho^\vee,b_-)$ for $2\rho^\vee=
\sum_{\al>0}\al^\vee.$
\medskip

The element $b_{-}= u_b(b)$ is a unique element
from $P_{-}$ that belongs to the orbit $W(b)$.
Thus the equality   $c_-=b_- $ means that $b,c$
belong to the same orbit. We will also use
$b_{+} \equal w_0(b_-),$ a unique element in $W(b)\cap P_{+}.$
In terms of $\pi_b,$
$$u_b\pi_b\ =\ b_-,\ \pi_b u_b\ =\ b_+.$$

Note that $l(\pi_b w)=l(\pi_b)+l(w)$ for all $b\in P,\ w\in W.$
For instance,
\begin{align}
&l(b_- w)=l(b_-)+l(w),\ l(wb_+)=l(b_+)+l(w),
\label{lupiw}
\\
&l(u_b\pi_b w)=l(u_b)+l(\pi_b)+l(w) \for b\in P,\,
 w\in W.\notag
\end{align}
\smallskip

{\bf Partial ordering on $P$.} 
It is necessary in the theory of nonsymmetric polynomials.
See \cite{O2,M4}. This ordering was also used in \cite{C2} in the
process of calculating the coefficients of $Y$\~operators.
The definition is as follows:
\begin{align}
&b \le c, c\ge b \for b, c\in P \iif c-b \in Q_+,
\label{order}
\\ &b \preceq c, c\succeq b \iif b_-< c_- \hbox{\ or\  }
\{b_-=c_- \hbox{\ and\ } b\le c\}.
\label{succ}
\end{align}
Recall that $b_-=c_- $ means that $b,c$
belong to the same $W$\~orbit.
We  write  $<,>,\prec, \succ$ respectively if $b \neq c$.

The following sets
\begin{align}
&\si(b)\equal \{c\in P, c\succeq b\},\
\si_*(b)\equal \{c\in P, c\succ b\}, \notag\\
&\si_-(b)\equal \si(b_-),\
\si_+(b)\equal \si_*(b_+)= \{c\in P, c_->b_-\}.
\label{cones}
\end{align}
are convex.
By {\em convex}, we mean that if
$ c, d= c+r\al\in \si$
for $\al\in R_+, r\in \Z_+$, then
\begin{align}
&\{c,\ c+\al,...,c+(r-1)\al,\ d\}\subset \si.
\label{convex}
\end{align}
\medskip

\subsection{More notations}
By  $m,$ we denote the least natural number
such that  $(P,P)=(1/m)\Z.$  Thus
$m=2 \for D_{2k},\ m=1 \for B_{2k} \and C_{k},$
otherwise $m=|\Pi|$.

We will need to include the case $t=0$ in our definition,
which requires minor deviations from the definitions
of \cite{C101},\cite{C4} and other author's papers.
Namely, we multiply all $T_i$ there by $t_i^{1/2}$
and change the formulas correspondingly.
 
The double affine Hecke algebra depends
on the parameters
$q, t_\nu,\, \nu\in \{\nu_\al\}.$ It will be defined
over the ring
$\Q[q^{\pm 1/m},t_\nu]$
formed by
polynomials in terms of $q^{\pm 1/m}$ and
$\{t_\nu \}.$ We will also use a greater ring
$$
\Q_{q,t}'\equal\{c\in \Q(q^{\pm 1/m},t_\nu)\, \mid\,
c \hbox{\ \, is well defined when\ \,} t_\nu=0\}, 
$$
which is a subring of the field of fractions of 
$\Q[q^{\pm 1/m},t_\nu]$

We set
\begin{align}
&   t_{\tal} = t_{\al}=t_{\nu_\al},\ t_i = t_{\al_i},\
q_{\tal}=q^{\nu_\al},\ q_i=q^{\nu_{\al_i}},\notag\\
&\where \tal=[\al,\nu_\al j] \in \tR,\ 0\le i\le n.
\label{taljx}
\end{align}

It will be convenient to use the parameters
$\{k_\nu\}$ together with  $\{t_\nu \},$ setting
$$
t_\al=t_\nu=q_\al^{k_\nu} \for \nu=\nu_\al, \and
\rho_k=(1/2)\sum_{\al>0} k_\al \al.
$$

Note that $(\rho_k,\al_i^\vee)=k_i=k_{\al_i}=
((\rho_k)^\vee,\al_i)$ for $i>0$; 
$(\rho_k)^\vee\equal \sum k_{\nu}(\rho_\nu)^\vee$.
Using that $w_0(\rho_k)=-\rho_k$, we obtain that
 $(\rho_k,-w_0(b))=(\rho_k,b)$. For instance,
$(\rho_k,b_+)=-(\rho_k,b_-)$, where 
$b_{+} \equal w_0(b_-)$ (see above).

By $q^{(\rho_k,\al)}$, we mean
$\prod_{\nu\in\nu_R}t_\nu^{((\rho_\nu)^\vee,\al)}$;
here $\al\in R$, $(\rho_\nu)^\vee=\rho_\nu/\nu$, and 
this product contains
only {\em integral} powers of $t_{\sht}$ and $t_{\lng}.$

For pairwise commutative $X_1,\ldots,X_n,$
\begin{align}
& X_{\tb}\ =\ \prod_{i=1}^nX_i^{l_i} q^{ j}
\iif \tb=[b,j],\ \hw(X_{\tb})\ =\ X_{\hw(\tb)}.
\label{Xdex}\\
&\hbox{where\ } b=\sum_{i=1}^n l_i \om_i\in P,\ j \in
\frac{1}{ m}\Z,\ \hw\in \hW.
\notag \end{align}
For instance, $X_0\equal X_{\al_0}=qX_\vth^{-1}$.

We set $(\tilde{b},\tilde{c})=(b,c)$ ignoring the affine extensions
in this pairing.
\medskip

\subsection{Main definition}
We note that $\pi_r^{-1}$ is $\pi_{r^*}$ and
$u_r^{-1}$ is $u_{r^*}$
for $r^*\in O\ ,$  $u_r=\pi_r^{-1}\om_r.$
The reflection $^*$ is
induced by an involution of the nonaffine Dynkin diagram
$\Gamma.$

\begin{definition}
The double affine Hecke algebra $\HH\ $
is generated over $\Q[q^{\pm 1/m},t_\nu]$ by
the elements $\{ T_i,\ 0\le i\le n\}$,
pairwise commutative $\{X_b, \ b\in P\}$ satisfying
(\ref{Xdex}),
and the group $\Pi,$ where the following relations are imposed:

(o)\ \  $ (T_i-t_i)(T_i+1)\ =\
0,\ 0\ \le\ i\ \le\ n$;

(i)\ \ \ $ T_iT_jT_i...\ =\ T_jT_iT_j...,\ m_{ij}$
factors on each side;

(ii)\ \   $ \pi_rT_i\pi_r^{-1}\ =\ T_j \iif
\pi_r(\al_i)=\al_j$;

(iii)\  $T_iX_b \ =\ X_b X_{\al_i}^{-1} \{t_iT_i^{-1}\} \iif
(b,\al^\vee_i)=1,\
0 \le i\le  n$;

(iv)\ $T_iX_b\ =\ X_b T_i$ if $(b,\al^\vee_i)=0
\for 0 \le i\le  n$;

(v)\ \ $\pi_rX_b \pi_r^{-1}\ =\ X_{\pi_r(b)}\ =\
X_{ u^{-1}_r(b)}
 q^{(\om_{r^*},b)},\  r\in O'$.
\label{double}
\end{definition}

Here and further the brackets $\{\cdot\}$ will be used
to {\em show explicitly} the elements from $t$\~localization
of $\HH$ that belong to $\HH\,$, i.e.,
those that do not involve $t_\nu^{-1}$ and other negative powers of 
$t_\nu$. It is not a new definition, but can help the readers to
see which operators are actually from $\HH\,$; 
quite a few (transitional) operators will involve
$t_\nu^{-1}$. We will postpone with the independent theory of
nil-DAHA, the limit of $\HH$ as $t\to 0$,
till the next paper(s). In this paper, we use the standard theory
of DAHA when convenient (which requires $t^{-1}$). The key 
examples are the elements $\{t_i^{1/2} T_i^{-1}\}$ which do belong
to $\HH$ thanks to the renormalization.  

One can rewrite (iii,iv) as in \cite{L}):
\begin{align}
&T_iX_b -X_{s_i(b)}T_i\ =\
(t_i-1)\frac{X_{s_i(b)}-X_b}
{X_{\al_i}-1},\ 0 \le i\le  n.
\label{tixi}
\end{align}

Given $\tw \in \tW, r\in O,\ $ the product
\begin{align}
&T_{\pi_r\tw}\equal \pi_r\prod_{k=1}^l T_{i_k},\where
\tw=\prod_{k=1}^l s_{i_k},
l=l(\tw),
\label{Twx}
\end{align}
does not depend on the choice of the reduced decomposition
(because $T_i$ satisfy the same ``braid'' relations
as $s_i$ do).
Moreover,
\begin{align}
&T_{\hv}T_{\hw}\ =\ T_{\hv\hw}\  \hbox{ whenever}\
 l(\hv\hw)=l(\hv)+l(\hw) \for
\hv,\hw \in \hW. \label{TTx}
\end{align}
In particular, we arrive at the pairwise
commutative elements: 
\begin{align}
& Y_{b} =  q^{(b_+-b,\,\rho_k)}
\prod_{i=1}^nY_i^{l_i} \iif
b=\sum_{i=1}^n l_i\om_i\in P,\ 
Y_i\equal T_{\om_i},b\in P.
\label{Ybx}
\end{align}
The factors here are needed to make them
from $\HH\,$;
$b_+$ is a unique element in $W(b)\cap P_+$.
Note that $Y_bY_{-b}= q^{2(b_+,\rho_k)}$.

Generally, if we replace $s_i$ by $T_i$
or $T_i^{-1}$ in any reduced decomposition
of $\hw\in \hW$, then such product belongs to $\HH\,$
upon the multiplication by
the product of $t_i$ corresponding to 
the terms $T_i^{-1}$. 

The relations dual to (iii,iv) hold
(for $i>0$ only):
\begin{align}\label{TYTL}
&\{t_i T^{-1}_i\}Y_b\ =\ Y_{s_i(b)} T_i\iif
(b,\al^\vee_i)=1,
\notag\\
& T_iY_b\ =\ Y_b T_i \iif (b,\al^\vee_i)=0,
 \ 1 \le i\le  n.
\end{align}
\smallskip

The counterpart of (\ref{tixi}) is as follows:
\begin{align}
&T_iY_b -Y_{s_i(b)}T_i\ =\
(t_i-1)\frac{Y_b-Y_{s_i(b)}}
{1-q^{-(\th',\rho_k)}Y_{-\al_i}},\ 1 \le i\le  n,
\label{tiyi}
\end{align}
where $\th'=\th,\vth$ respectively for long, short $\al_i$
(it is the only root in the intersection $W(\al_i)\cap P_+$).

Here and below we use that given $b\in P$,
replacing all $T_i^{\pm 1}$ by $t_i^{\pm 1}$ in the product of 
(\ref{Ybx}) for $Y_b$ results in the $t$\~power
$q^{2(\rho_k,b)}=\prod_\nu t_\nu^{2((\rho_\nu)^\vee,\,b)}$.

In the standard DAHA theory,
$q^{-(b_+,\,\rho_k)}Y_b$ for any $b$ can be represented 
as the product $\pi_r (t_l^{\mp 1/2}T_{i_l}^{\pm 1})\cdots 
(t_1^{\mp 1/2}T_{i_1}^{\pm 1})$ for
a given reduced decomposition $b=\pi_r s_{i_l}\cdots s_{i_1}$
and proper choice of $\{\pm\}$. The number of terms
is $l=l(b)=2(\rho^\vee,\,b_+)$. Only positive
powers $T_i^{+1}$ will appear in this product when $b\in P_+$. 
The total number of the terms $T_i^{\pm 1}$ with 
$\nu_i=\nu$ in this product 
equals $2((\rho_\nu)^\vee,\,b_+)$.
\medskip

\setcounter{equation}{0}
\section{Polynomial representation}
From now on, we will switch from $\HH\ $ to its
{\em intermediate subalgebra} $\HH^\flat\subset\HH\ $
with $P$ replaced by a lattice $B$ between $Q$ and $P$
(see \cite{C12}). Accordingly, $\Pi$ is
changed to the preimage $\Pi^\flat$ of $B/Q$ in $\Pi.$ Generally,
there can be two different
lattices $B_X$ and $B_Y$ for $X$ and $Y.$
We consider only $B_X=B=B_Y$ in the paper;
respectively, $a,b\in B$ in $X_a,Y_b.$

We also set $\hW^{\flat}=B\cdot W\subset \hW,\ $
and replace $m$ by the least $\tilde{m}\in \N$ such that
$\tilde{m}(B,B)\subset \Z$ in the definition of the
$\Q_{q,t}'.$

Note that $\HH^\flat$ and the polynomial representations
(and their rational and trigonometric
degenerations)  
are actually defined over $\Z$ extended by the parameters of DAHA.
However the ring $\Q_{q,t}'$ will be sufficient in this paper.
\smallskip

The {\em Demazure-Lusztig operators}
are as follows:
\begin{align}
&T_i\  = \  t_i s_i\ +\
(t_i-1)(X_{\al_i}-1)^{-1}(s_i-1),
\ 0\le i\le n;
\label{Demazx}
\end{align}
they obviously preserve $\Q[q,t_\nu][X_b]$.
We note that only the formula for $T_0$ involves $q$:
\begin{align}
&T_0\  = \ t_0s_0\ +\ (t_0-1)
(X_0 -1)^{-1}(s_0-1),\hbox{\ where\ }\notag\\
&X_0=qX_\vth^{-1},\
s_0(X_b)\ =\ X_bX_{\vth}^{-(b,\vth)}
 q^{(b,\vth)},\
\al_0=[-\vth,1].
\end{align}

The map sending $ T_j$ to the corresponding operator from
(\ref{Demazx}), $X_b$ to $X_b$
(see (\ref{Xdex})) and 
$\pi_r\mapsto \pi_r$ induces a
$ \Q_{ q,t}'$\~linear
homomorphism from $\HH^\flat\, $ to the algebra of
linear endomorphisms
of $\Q_{ q,t}'[X]$.
This $\HH^\flat\,$-module is faithful
and remains faithful when  $q,t$ take
any complex values assuming that
$q\neq 0$ is not a root of unity.
It will be called the
{\em polynomial representation};
the notation is
$$
\v\equal \Q_{q,t}'[X_b]\ =\ \Q_{q,t}'[X_b, b\in B].
$$

The images of the $Y_b$ are called the
{\em difference-trigonometric Dunkl operators}.

The polynomial representation
is the $\HH^\flat\,$\~module induced from the one-dimensional
representation $T_i\mapsto t_i,\,$ 
$Y_b\mapsto q^{2(\rho_k,b)}$
of the affine Hecke subalgebra 
$\h_Y^\flat=\lan T_i,Y_b\ran.$ Here we extend the ring of
constants to $\Q_{q,t}'$.
\medskip

\subsection{Macdonald polynomials}
There are two equivalent definitions of the
{\em nonsymmetric Macdonald polynomials\,}, denoted by
$E_b(X) = E_b^{(k)}$ for $b\in B$; they belong to
$\Q(q,t)'[X_a,a\in B]$. The first definition is based
on the truncated theta function due to Macdonald:
\begin{align}
&\mu(X;t) = \mu^{(k)}(X)=\prod_{\al \in R_+}
\prod_{j=0}^\infty \frac{(1-X_\al q_\al^{j})
(1-X_\al^{-1}q_\al^{j+1})
}{
(1-X_\al t_\al q_\al^{j})
(1-X_\al^{-1}t_\al^{}q_\al^{j+1})}.\
\label{mu}
\end{align}

We will mainly consider $\mu$ as a Laurent series with the
coefficients in  the ring $\Q[t_\nu][[q_\nu]]$ for
$\nu\in \nu_R=\{\nu_{\sht},\nu_{\lng}\}$.
The constant term of a Laurent series $f(X)$ will be denoted
by $\langle  f \rangle.$ Then
\begin{align}
&\langle\mu\rangle\ =\ \prod_{\al \in R_+}
\prod_{j=1}^{\infty} \frac{ (1- q^{(\rho_k,\al)+j\,\nu_\al})^2
}{
(1-t_\al q^{(\rho_k,\al)+j\,\nu_\al})
(1-t_\al^{-1}q^{(\rho_k,\al)+j\,\nu_\al})
}.
\label{consterm}
\end{align}
Recall that $q^{(z,\al)}=q_\al^{(z,\al^\vee)},
\ t_\al=q_\al^{k_\al}.$
This equality is equivalent
to the Macdonald constant term conjecture
proved in complete generality in \cite{C2}.

Let $\mu_\circ\equal \mu/\langle \mu \rangle$.
The coefficients of the Laurent series $\mu_\circ$ are
from the ring $\Q_{q,t}'$.
\medskip

The polynomials $E_b$ are uniquely determined from
the relations 
\begin{align}
&E_b-X_b\ \in\ \oplus_{c\succ b}\Q_{q,t}' X_c,\
\lan E_b X_{c}^{-1}\mu_\circ\ran = 0 \for B \ni c\succ b.
\label{macd}
\end{align}
for generic $q,t$ and form a basis in $\Q_{q,t}'[X_b]$.

This definition is due to Macdonald (for
$k_{\sht}=k_{\lng}\in \Z_+ $), who extended
the construction from \cite{O2}.
The general (reduced) case was considered in \cite{C4}.

Another approach is based on the $Y$\~operators.
We continue using the same notation $X,Y,T$
for these operators acting in the polynomial
representation. Let $X_a(q^{b})\ =\
q^{(a,b)}$ for  $a,b\in P.$

\begin{proposition}
The polynomials $\{E_b, b\in B\}$
are unique (up to proportionality) eigenfunctions of
the operators $Y_a$ ($a\in P$) 
acting in $\Q_{q,t}'[X]:$
\begin{align}
&Y_a(E_b)\ =\ q^{(a_+,\,\rho_k)-(a,b_\#)}\,E_b\, \hbox{\ for\ }\,
b_\#\equal b- u_b^{-1}(\rho_k),
\label{Yone} 
\end{align}
$u_b=\pi_b^{-1}b$ is from Proposition \ref{PIOM},\
$b_\#=\pi_b(\!(-\rho_k)\!)$.
\label{YONE}
\end{proposition}
\qed

The coefficients of the Macdonald
polynomials are well known to be
rational functions in terms of $q_\nu,t_\nu$. In this
paper, we use that these coefficients actually belong 
to $\Q_{q,t}'$, i.e., are
well defined when $t_\nu=0$ for all $\nu$. It readily follows
from (\ref{epolexists}) below.

\medskip
\subsection{Symmetric polynomials}
Following Proposition \ref{YONE}, the
{\em symmetric Macdonald polynomials\,} $\ P_b=P_b^{(k)}\ $
can be introduced as
eigenfunctions of the $W$\~invariant
{\em difference} operators 
\begin{align}
&L_{a_+}=
\hbox{Red}_W(\sum_{a'\in W(a_+)}Y_{a'}) \for a_+\in B_+\,,
\label{maqsum}
\end{align}
where $\hbox{Red}_W$ is the
restriction to the space  $\v^{\,W}$ of $W$\~invariants
of $\v$. Explicitly,
\begin{align}
L_{a_+}&(P_{b_-})=q^{(a_+\,,\, \rho_k)}
(\sum_{a'\in W(a_+)}q^{-(a',\ b_--\rho_k)})\, 
P_{b_-},\ 
b_-\in B_-,\notag\\
&P_{b_-}=\sum_{b\in W(b_-)}X_b \mod
\oplus_{c_-\succ b_-}\Q(q,t) X_c.\
\label{Lf}
\end{align}

These polynomials were introduced in \cite{M2,M3}.
They were used for the first
time in Kadell's unpublished work (classical root systems).
In the case of $A_1$, they are due to Rogers.

The connection between $E$ and $P$ is as follows
\begin{align}
&P_{b_-}\ =\ \P_{b_+} E_{b_+}, \ b_-\in B_-,\
b_+=w_0(b_-),
\notag\\
&\P_{b_+}\equal\sum_{c\in W(b_+)} T_{w_c}, \where
\label{symmetr}
\end{align}
$w_c\in W$ is the element of the least length such that
$c=w_c(b_+)$. Taking the complete $t$\~symmetrization
$\P$ here (with the summation over all $w$), one obtains
$P_{b_-}$ up to proportionality. See \cite{O2,M4,C4}. 
\medskip

There are two different kinds of inner products in 
$\v$ from \cite{C101} and other works. In the symmetric
setting, they essentially coincide. We will need here only 
the inner products of the symmetric polynomials
$P_b$ for $b=b_-$\,:
\begin{align}\label{normppols}
\lan P_{b}(X)P_{c}(X^{-1})\mu_\circ\ran&\\
=\ \de_{bc}\prod_{\al>0}\prod_{j=0}^{-(\al^\vee,b)-1}&
\Bigl(\frac{
(1-q_\al^{j+1}t_\al^{-1} X_\al(q^{\rho_k}))
(1-q_\al^j t_\al  X_\al(q^{\rho_k}))
}{
(1-q_\al^j X_\al(q^{\rho_k})) (1-q_\al^{j+1} X_\al(q^{\rho_k}))
}\Bigr).\notag
\end{align}
\medskip

\subsection{Using intertwiners}
The following map can be uniquely extended to
an automorphism of $\HH^\flat\,$ where 
proper fractional powers of $q$ are added
(see \cite{C15},\cite{C4},\cite{C12}):
\begin{align}
& \tau_+:\  X_b \mapsto X_b, \
\pi_r \mapsto q^{-\frac{(\om_r,\om_r)}{2}}X_r\pi_r,
\ Y_r \mapsto X_rY_r q^{-\frac{(\om_r,\om_r)}{2}},
\notag\\
& \tau_+:\ T_0\mapsto X_0^{-1} \{t_0T_0^{-1}\},\ 
Y_\vth \mapsto  X_0^{-1} \{t_0 T_0^{-1}\}T_{s_\vth}. 
\label{taux}
\end{align}
This automorphism  fixes
$T_i\,(i\ge 1)$,$\ t_\nu,\ q$
and fractional powers of $q.$ 
\medskip

The $Y$\~intertwiners serve as creation operators
in the theory of nonsymmetric Macdonald polynomials.
Following \cite{C1,C101}, let
\begin{align}
&\Psi_i^c=
\tau_+(T_i) + (t_i-1)
(X_{\al_i}(q^{c_{\#}})-1)^{-1}, 0\le i\le n.
\label{Phijb}
\end{align}
We will use the pairing from (\ref{dform})
and the affine action $\hw\llb c\rrb$
from (\ref{afaction}).

\begin{theorem}\label{PHIEB}
Given $c\in B,\ 0\le i\le n$ such that
$(\al_i, c+d)> 0,$
\begin{align}
&q^{(c,c)/2-(b,b)/2} E_{b} \ =\ \Psi_i^c(E_c)
\for b= s_i\llb c\rrb.
\label{Phieb}
\end{align}
If  $(\al_i, c+d)=0,$ then
\begin{align}
&\tau_+(T_i) (E_c) \ =\ t_i E_c, \ 0\le i\le n,
\label{Tjeco}
\end{align}
which  results in the relations $s_i(E_c)=E_c$ as
$i>0$.  For $b=\pi_r\llb c\rrb,$ where
the indices $\,r\,$ are from $O',$
\begin{align}
&q^{(c,c)/2-(b,b)/2}E_b\ =\ \tau_+(\pi_r)(E_c)\ =\
X_{\om_r}q^{-(\om_r,\om_r)/2}\pi_r(E_c).
\label{pireb}
\end{align}
Also $\tau_+(\pi_r)(E_c)\neq E_c$ for $\pi_r\neq$id,
since $\pi_r\llb c\rrb\neq c$ for any $c\in B$.

\sq
\end{theorem}
\smallskip

If  $(\al_i, c)>0$ and $i>0$, then
the set $\la(\pi_b)$ is obtained from
$\la(\pi_c)$ by adding $[\al,(c_-,\al)]$
for $\al=u_c(\al_i)\in R_-$
and $(c_-,\al^\vee)=(c,\al_i^\vee)>0.$
When $i=0$ and $(\al_0, c+d)=-(c,\vth)+1>0$,
then the root $[\al,(c_-,\al)+1]$ is added to 
$\la(\pi_c)$ for $\al= u_c(-\vth)=\al^\vee\in R_-$
and $( c_-, \al)=-(c,\vth)\ge 0.$

In each of these two cases,
$(\al_i,u_c^{-1}(\rho))=(\al,\rho)<0$ and
the powers of $t_\nu$ in  
\begin{align}\label{xaliq}
&X_{\al_i}(q^{c_{\#}})=q^{(\al_i,c-u_c^{-1}(\rho_k)+d)}=
q^{(\al_i,c+d)}\prod_\nu t_\nu^{-(\al,(\rho_\nu)^\vee)}
\end{align}
are from $\Z_+$ with that of $t_{i}$ strictly 
positive. 
\medskip

Due to Theorem \ref{PHIEB}
(see also \cite{C1}, Corollary 5.3),
the polynomial $E_b$ exists if
\begin{align}\label{epolexists}
&\prod_{[-\al,\,\nu_\al\,j]\in \la\,(\pi_b)}
\bigl(
1- q_\al^{j}X_\al(q^{\rho_k})\bigr)\neq 0.
\end{align}
If $b\in B_-$ and the latter inequality holds for
$b_+=w_0(b)\in B_+,$
then the symmetric polynomial
$P_b$ is well defined. If $t_\nu=0$ for all $\nu$,
then $E$\~polynomials and $P$\~polynomials are
always well defined,
which gives that their coefficients are polynomials
in terms of $q$. 

\subsection{Spherical polynomials}
The following renormalization
of the $E$-polynomials is of major importance in
the Fourier analysis (see \cite{C4}):
\begin{align}\label{ebebs}
\e_b\ \equal&\ E_b(X)(E_b(q^{-\rho_k}))^{-1},\where b\in B,
\\
E_{b}(q^{-\rho_k}) \ =&\ q^{(\rho_k,b_-)}
\prod_{[\al,j]\in \la\,'\,(\pi_b)}
\Bigl(
\frac{
1- q_\al^{j}t_\al X_\al(q^{\rho_k})
 }{
1- q_\al^{j}X_\al(q^{\rho_k})
}
\Bigr).\notag
\end{align}
This definition requires the $t$\~localization.

We call them nonsymmetric
{\em spherical polynomials\,}.
Formula (\ref{ebebs}) is the Macdonald
{\em evaluation conjecture}
in the nonsymmetric variant from \cite{C4}.
See \cite{C3} for the symmetric evaluation conjecture.

The following {\em duality formula} holds
for $b,c\in B\, :$
\begin{align}
&\e_b(q^{c_{\#}})\ =\ \e_c(q^{b_{\#}}),\
b_\# = b- u_b^{-1}(\rho_k),
\label{ebdual}
\end{align}
which is the main justification of the definition of $\e_b$.

Given $b\in B,$
the polynomial
$\e_b$
is well defined for $q,t\in \C^*$ if
\begin{align}
&\prod_{[\al,j]\in \la\,'\,(\pi_b)}
\bigl(
1- q_\al^{j}t_\al X_\al(q^{\rho_k})\bigr)\ \neq\ 0.
\label{esphexists}
\end{align}
\smallskip

In the symmetric setting,
\begin{align}\label{pebebs}
\p_b\ \equal&\ P_b(X)(P_b(q^{-\rho_k}))^{-1} \where b\in B_-\,,
\\
P_{b}(q^{-\rho_k})=P_{b}(q^{\rho_k}) \ =&\ q^{(\rho_k,b_-)}
\prod_{\al>0}\prod_{j=0}^{-(\al^\vee,b)-1}
\Bigl(
\frac{
1- q_\al^{j}t_\al X_\al(q^{\rho_k})
 }{
1- q_\al^{j}X_\al(q^{\rho_k})
}
\Bigr).\notag
\end{align}
The symmetric duality reads as follows:
\begin{align}
&\p_b(q^{c-\rho_k})\ =\ \p_c(q^{b-\rho_k}),\
\for b,c\in B_-\,.
\label{pebdual}
\end{align}
\medskip

The norm formula becomes entirely conceptual:
\begin{align}\label{normpipols}
&(\lan \p_{b}(X)\p_{b}(X^{-1})\mu_\circ\ran)^{-1}\ 
=\ \sum_{a\in W(b)}\,\mu(\pi_{a})\mu(\hbox{id})^{-1},\\
& \where \mu(\hw)\equal\mu(\hw(\!(q^{-\rho_k})\!))
\for \hw\in \hW.\notag
\end{align}
It is a direct corollary of the fact that the Fourier
transform sends the $\p$\~polynomials to the delta-functions;
see \cite{C101}.
\medskip

\subsection{The limit 
\texorpdfstring{\mathversion{bold}$t\to 0$}{\em t=0}}
Let $\overline{\HH}^{\,\flat}\,$ by the reduction of
$\HH^\flat\,$ for $t_\nu=0$, where $\nu\in \nu_R$. 
It can be called the {\em nil-DAHA} or
the {\em crystal DAHA}. The polynomials
$\overline{E}_b,\overline{P}_{b_-}$ are well defined
and linearly generate 
$\overline{\v}$ and $\overline{\v}\,^W$
correspondingly;  $\overline{\v}=\Q_q[X_b,b\in B]$,
where $\Q_q=\Q(q^{1/m})$ is the ring of {\em polynomials}
in terms of $q^{1/m}$ with $m$ from the definition of
$\Q_{q,t}'$.
We will denote $T_i(t=0)$ by $\overline{T}_i$.

Theorem \ref{PHIEB} holds under this specialization
and gives quite a constructive approach to the
$\overline{E}$\~polynomials. The intertwiners 
$\Psi_i^c$ from (\ref{Phijb})
that appear in the formulas for $\overline{E}_b$ are all
in the form $\tau_+(\overline{T}_i)+1$ in this limit. 
It is directly connected with the fact that
$T_i'=\overline{T}_i+1$ satisfy the same homogeneous
Coxeter relations as $\{T_i,\, 0\le i\le n\}$ do,  
a special feature of the nil-DAHA. It readily results from the
theory of intertwiners, but, of course, can be checked
directly too.

The action 
of $\pi_r$ on $\{T_i'\}$ by conjugation
obviously remains unchanged. Thus relations 
(i,ii) from Definition \ref{double} hold and,
given $\hw\in \hW$, 
the element $T'_{\hw}=\pi_rT'_{i_l}\cdots T'_{i_1}$
does not depend on the choice of the reduced
decomposition $\hw=\pi_r s_{i_l}\cdots s_{i_1}$.
For instance, operators $\Pi'_i\equal\tau_+(T'_{-\om_i})$ for
$i=1,\ldots,n$ are pairwise commutative and, importantly,
$W$\~invariant. 

Indeed, one has\,:\,
$\Pi'_{b}=\prod_{i=1}^n\, (\Pi'_i)^{n_i}$ for
$B_-\ni b=-\sum n_i\,\om_i$. Provided that all $n_i>0$,
the {\em reduced} decomposition $b=b_-=w_0\pi_{b_+}$ holds
for the longest element $w_0\in W$ and $b_+=w_0(b)\in B_+$;
see (\ref{lupiw}).
Thus $\Pi'_b$ is divisible on the left by $(\overline{T}_i+1)$
for any $i>0$ and therefore divisible by the $W$\~symmetrizer on the
left. It results in the $W$\~invariance of $P_b$ for any
$b\in B_-$.   

The $W$\~invariance of $\{\Pi_b,\,b\in B_-\}$ simplifies
significantly the relation of the $\overline{E}$\~polynomials to
the $\overline{P}$\~polynomials:
\begin{align}\label{eplimbar}
&\overline{P}_{b}\ =\ \overline{E}_{b} \mbox{\ \ for\ \ }
b=b_-\in B_-\,.
\end{align}
In more detail, we have the following explicit
proposition.

\begin{proposition} \label{BAREP}
(i) In the representation $\overline{\v}$
of $\overline{\HH}^{\,\flat}\,$,
the polynomial\,  $\tau_+(T'_{\hw})(1)$ 
equals \,  $q^{r_b}\,\overline{E}_b$ \,for $\hw=\pi_b,\, b\in B$, 
\, $r_b\in \Q$.

(ii) In the symmetric case,
\begin{align}\label{barpform}
&\Pi'_b(1)\ =\ q^{r_b}\overline{P}_b \for b\in B_-\,,\ r_b\in\Q,
\end{align} 
where $\Pi'_i$ can be replaced by their restrictions
$\hbox{Red}_{\,W}(\Pi'_i)$ to $\overline{\v}^{\,W}$, 
which are pairwise commutative $W$\~invariant
difference operators.\sq
\end{proposition}
\smallskip

It is important that only positive powers of $q$ appear in the 
coefficients of $\overline{E}_b$. The coefficients of these
$q$\~polynomials are non-negative integers. One can obtain
it from the interpretation via Demazure characters or 
using the intertwiners (we are going to discuss it in further
papers). As $q\to 0$, the polynomials 
$\overline{P}_{b_-}$ become
the classical finite dimensional Lie characters,
which can be seen, for instance, from (\ref{normppolsbar})
below.

For the affine root systems considered in this paper
(with $\al_0$ in terms of the maximal {\em short} root $\vth$),
the connection was established
between the polynomials $E_b(t\to\infty)$ and
the Demazure characters of the corresponding
irreducible affine Lie algebras. See 
\cite{San} and, especially, \cite{Ion1}, Theorem 1.
Paper \cite{Ion1} is based on the technique of 
intertwiners (from \cite{KnS} in the $GL_n$\~case
and \cite{C1} for arbitrary reduced root systems). 
We will not discuss this important relation in this paper.

There is a relation between the limit $t\to 0$
used here and the one $t\to \infty$. It 
goes through the general formula
\begin{align}
&E_b^*\ =\
\prod_{\nu\in \nu_R} t_\nu ^{l_\nu(u_b)-l_\nu(w_0)}\,
T_{w_0}(E_{\varsigma(b)}), \hbox{\ \ where}\label{ebast}\\
&X^*=X^{-1}, q^*=q^{-1}, t^*=t^{-1},\ 
\varsigma(b)=-w_0(b),\notag
\end{align} 
form \cite{C101} and other author's works.
This connection is especially
simple for the symmetric polynomials: 
$P_{b}(X)^*=P_{b}(X^{-1})$ for $b=b_-$, i.e.,
$\overline{P}_{b}=P_{b}(t\to 0)=P_{\varsigma(b)}(t\to \infty).$
We use that $P_b(X^{-1})=P_{\varsigma(b)}(X)$. This connects
the $q$\~Hermite polynomials and the Demazure characters 
for $b=b_-$.
\smallskip

Concerning the orthogonality of
$\overline{P}$, 
the denominator of the $\mu$\~function from
(\ref{mubar}) vanishes in the
limit:  
\begin{align}
&\overline{\mu}\ = \ \prod_{\al \in R_+}
\prod_{j=0}^\infty (1-X_\al q_\al^{j})
(1-X_\al^{-1}q_\al^{j+1}).
\label{mubar}
\end{align}
The constant term formula becomes a well-known
identity:
\begin{align}
&\langle\overline{\mu}\rangle\ =\ \prod_{i=1}^{n}
\prod_{j=1}^{\infty} \frac{1}
{1-q_i^j}\,, \where q_i=q^{\nu_i}.
\label{constermbar}
\end{align}

For $b,c\in B_-\,$,
the norm formula from (\ref{normppols}) reads as:
\begin{align}\label{normppolsbar}
&\lan \overline{P}_{b}(X)\overline{P}_{c}(X^{-1})
\overline{\mu}_\circ\ran\ 
=\ \de_{bc}\prod_{i=1}^n\prod_{j=1}^{-(\al_i^\vee,b)}
(1-q_i^{j})\,.
\end{align} 
\medskip

\setcounter{equation}{0}
\section{Spherical and Whittaker functions}
We will begin with the identities involving
the Gaussians, which are essentially from \cite{C5}; 
then their limits $t\to 0$ will be considered.

The second part of this section is devoted to
the Whittaker limit of the
$q,t$\~spherical function from \cite{C5}, which results
in a formula for the $q$\~Whittaker function
in terms of the $\overline{P}$\~polynomials.  

We note that the {\em Whittaker limit} is a general 
procedure that can be applied to any solutions
of the Macdonald eigenvalue problem (and its
various degenerations and generalizations).
\medskip

\subsection{Gauss-type integrals}
By the {\it Gaussians} $\tga$ we mean 
\begin{align}
&\tga^{\oplus\,}= \sum_{b\in B} q^{-(b,b)/2}X_b,\ 
\tga^{\ominus\,}= \sum_{b\in B} q^{(b,b)/2}X_b.
\label{gauser}
\end{align}
The multiplication by $\tga^{\ominus\,}$ preserves 
the space of Laurent series
with coefficients in $\Q[t][[q^{\frac{1}{2\tilde{m}}}]]$, where
$\tilde{m}(B,B)= \Z$ is from the definition of 
$\Q_{q,t}'.$ Accordingly, the coefficients must be
taken from  $\Q[t][[q^{-\frac{1}{2\tilde{m}}}]]$ when
the Gaussian $\tga^{\oplus}$ is taken.

We will also use the {\em real Gaussians} defined as
\begin{align}\label{gaussreal}
&\ga^{\pm 1}=q^{\pm x^2/2}, \where
X_b\equal q^{x_b}, x_b=(x,b), 
x^2=\sum_i\, x_{\al_i}\,x_{\om_i^\vee}.
\end{align}

Note that considering $\tga^{\oplus,\ominus}$ 
as holomorphic functions (provided that
$|q|>1$ and, respectively, $|q|<1$) the functions
$\tga^{\oplus\,}/\ga$ and $\tga^{\ominus\,}\ga$ are 
$B$\~periodic in terms of $x$.
\medskip

The {\em $q$\~Mehta\~Macdonald identity} from \cite{C5}
\begin{align}
&\langle \tga^{\ominus\,}\mu_\circ\rangle\ =\
\prod_{\al\in R_+}\prod_{ j=1}^{\infty}\Bigl(\frac{
1-t_\al^{-1} q_\al^{(\rho_k,\al^\vee)+j}}{
1-           q_\al^{(\rho_k,\al^\vee)+j} }\Bigr)
\label{mehtamu}
\end{align}
provides the normalization constant
for the $q$\~Gauss integrals 
\begin{align}\label{pbgauss}
\langle P_b(X) P_c(X) \tga^{\ominus\,}\mu_\circ \rangle& \\  
=\ q^{\frac{(b,b)+(c,c)}{2} -(b+c\,,\,\rho_k)}& 
P_c(q^{b-\rho_k})P_b(q^{-\rho_k})
\langle \tga^{\ominus\,}\mu_\circ\rangle,\notag
\end{align} 
where $b,c\in B_-\,$. Obviously, it implies the
duality formula  (\ref{pebdual}).
Formula (\ref{pbgauss}) can be naturally extended
to the $E$\~polynomials (the proof even
becomes simpler), but we do not need it in this paper.

There are counterparts of (\ref{pbgauss}) for $\tga^{\,\oplus}$
(treated as an analytic function for $|q|>1$),
and for the Jackson summation taken instead of the constant term
functional. See \cite{C5,C101}.
The considerations from this paper
can be readily extended to these cases.
\medskip

{\bf Taking the limit.}
Let us tend $t\to 0$ in (\ref{pbgauss}). The
definition of the $P$\~polynomials implies that
\begin{align}\label{limpct}
&\lim_{t\to 0} q^{-(c,\,\rho_k)} 
P_c(q^{z-\rho_k})\ =\ q^{(c_+\,,\,z)} \for c\in B_-\,, c_+=w_o(c).
\end{align}
Note that $-(c,\rho_k)=(c_+,\rho_k)$.
For instance, it matches the evaluation
formula in (\ref{pebebs}): 
$\lim_{t\to 0} q^{-(c,\rho_k)} 
P_c(q^{-\rho_k})=1$.

We come to the following formulas ($c\in B_-\,$):
\begin{align}
&\langle \tga^{\ominus\,}\overline{\mu}_\circ\rangle\ =\
\prod_{i=1}^n\prod_{ j=1}^{\infty}(1-q_i^j),
\label{mehtamul}
\end{align} 
\begin{align}\label{pbgaussl}
&\langle \overline{P}_b(X) \overline{P}_c(X) 
\tga^{\ominus\,}\overline{\mu}_\circ \rangle\  =\  
 q^{\frac{(b,b)+(c,c)}{2}}\, X_{c_+}(q^{b})\,
\langle \tga^{\ominus\,}\overline{\mu}_\circ\rangle.
\end{align} 
Here $X_{c_+}(q^{b})=q^{(c_+,b)}=q^{(c,b_+)}=X_{b_+}(q^{c})$.
\medskip

\subsection{Global spherical function}
One of the  main advantages of the technique of Gaussians is
a possibility to introduce the spherical function as
a reproducing kernel of the Fourier transform from
$\v\ga^{-1}$, the polynomial representation multiplied
by the Gaussian $\ga^{-1}$, 
to the $\HH^\flat\,$\~module  $\v\ga$. We
will need only the symmetric case here. We assume
that  $|q|<1$, which makes the
considerations ``naturally" compatible with the limit $t\to 0$.
In this setting, the construction below is directly 
related to the identities  (\ref{pbgauss})
(correspondingly, (\ref{pbgaussl}) in the limit).

We note that if the whole polynomial representation is considered,
then the corresponding anti-involutions of $\HH^\flat$,
generally, require the $t$\~localizations. Correspondingly, 
the definition of the Whittaker limit of 
the nonsymmetric counterpart of formula (\ref{pexla}) 
below (see \cite{C5}) becomes more subtle.

We will use the notation $\tga_\la$ and $\ga_\la$
for the Gaussians defined for another set of variables
$\La$ completely analogous to $X$ ($\tga_x,\ga_x$ are old 
$\tga,\ga$). Thus,  $\tga_\la=\tga(q^\la)$ and 
$\ga_\la=\ga(q^\la)=q^{\la^2/2}$.
We will also use
\begin{align}
&\langle \ga\rangle_{\rho_k}\ \equal\  
\sum_{a\in B} q^{\frac{(\rho_k+a,\rho_k+a)}{2}}=
\tga^{\ominus}(q^{\rho_k})q^{\frac{(\rho_k,\rho_k)}{2}}
\label{lgarhok}. 
\end{align}

\begin{theorem}\label{GLOBSPH}
Provided that $|q|<1$, the function $\Psi$ from the relation
\begin{align}
&\tga_x^{\ominus}\tga_\la^{\ominus}
\mathfrak{P}^\circ(X,\La)/\tga_x^{\ominus}(q^{\rho_k})
\notag\\
=\ \Psi(X,\La;q,t)&\equal
\sum_{b\in B_-} q^{\frac{(b,b)}{2} -(\rho_k,b)} 
 \frac{P_b(X)\ P_b(\La^{-1})} {
\langle P_b(X) P_b(X^{-1})\mu_\circ\rangle}
\label{pexla}
\end{align}
is a well-defined Laurent series.
It is an analytic function for all $X,\La$ and for
any choice of $t_\nu$ assuming that all $P$\~polynomials
exist (the conditions $|t_\nu|<1$ are sufficient).

The function $\mathfrak{P}^\circ(X,\La)$ defined
via (\ref{pexla}) is meromorphic
for all $X,\La$ and analytic 
apart from the zeros of $\tga_x^{\ominus\,}\tga_\la^{\ominus\,}$.
Replacing $\tga_x^{\ominus\,}\tga_\la^{\ominus\,}$ by 
$\ga_x^{-1}\ga_\la^{-1}$ in this definition, the corresponding
function will be denoted simply by $\mathfrak{P}(X,\La)$;
it becomes totally analytic but not a (single-valued)
function in terms of $X_b,\La_b$.

Both functions, $\mathfrak{P}^\circ(X,\La)$ and 
$\mathfrak{P}(X,\La)$,  are
 $X\leftrightarrow \La$\~symmetric, $W$\~invariant
with respect to $X$ and $\La$ and satisfy 
the following extension of the eigenvalue problem 
from (\ref{Lf}):
\begin{align}
L_{a_+}&(\mathfrak{P}(X,\La))\ =\ q^{(a_+\,,\, \rho_k)}
(\,\sum_{a'\in W(a_+)}\,\La_{a'}^{-1})\, 
\mathfrak{P}(X,\La).
\label{LfLa}
\end{align}\sq
\end{theorem} 

We note that $P_b(X^{-1}) P_b(\La)=
P_{\varsigma{b}}(X) P_{\varsigma{b}}(\La^{-1})$ in 
(\ref{pexla}); recall that $\varsigma(x)
=-w_0(x)$ and $P_{\varsigma{b}}(X)=P_b(X^{-1})$.
Applying $\varsigma$ to the summation
index $b$ does not change the result. Thus:
$$
\mathfrak{P}^\circ(X,\La)=\mathfrak{P}^\circ(\La,X)=
\mathfrak{P}^\circ(\varsigma(X),\varsigma(\La)).
$$ 

The following can be used for an {\em abstract} (i.e., without
an explicit formula) definition of the 
function $\mathfrak{P}^\circ(X,\La)$.
It goes through
the spherical polynomials 
$\{\p_c=P_c/P_c(q^{-\rho_k}), c\in B_-\}$ with a common
coefficient of proportionality:
\begin{align}
& \mathfrak{P}^\circ(X,q^{c-\rho_k})
\ =\  
\frac{P_c(X)}{P_c(q^{-\rho_k})} 
\prod_{\al\in R_+}\prod_{ j=1}^{\infty}\Bigl(\frac{ 
1- q_\al^{(\rho_k,\al^\vee)+j}}{
1-t_\al^{-1}q_\al^{(\rho_k,\al^\vee)+j} }\Bigr). 
\label{hatmuxsym}
\end{align}
Here we substitute $\la=c_\#,\, \La=q^{c_\#}$\, in the left-hand
side of (\ref{LfLa}) and divide it by the Gaussian
$\tga^{\ominus}_x$. This formula can be considered
as a $q,t$\~generalization of the Shintani-Casselman-Shalika
formula from \cite{Shi,CS}. Its limit as $t\to 0$ will be
discussed in the next section.
\medskip

\subsection{Global Whittaker function}
We are now in a position to define the {\em global\, 
$q$\~Whittaker function} 
$\widetilde{\mathfrak{P}}^\circ_x(X,\La)$
from the relation
\begin{align}\label{whitgen}
& \tga_x^{\ominus}\widetilde{\mathfrak{P}}^\circ_x(X,\La)\ \equal\ 
\lim_{t\to 0} \frac{\tga^{\ominus}(q^{x-\rho_k})}
{\tga^\ominus(q^{\rho_k})}
\ \mathfrak{P}^\circ (q^{-\rho_k}X,\La).
\end{align}

Here we always assume that $t_\nu\to 0$ for all $\nu$.
The function $\widetilde{\mathfrak{P}}_x$ is
defined for $\ga^{-1}$ instead of $\tga^{\ominus}$:
\begin{align}\label{whitgenmin}
& \widetilde{\mathfrak{P}}_x(X,\La)\ \equal\ 
\lim_{t\to 0} q^{(x\,,\,\rho_k)}
{\mathfrak{P}}(q^{-\rho_k}X,\La).
\end{align}
More explicitly, provided that $|q|<1$ (we will not show the
dependence of $\La$ here and where it cannot lead to  
misunderstanding):
\begin{align}\label{whitgenn}
& \widetilde{\mathfrak{P}}_x(X)\ \equal\ 
\ga_x\, \lim_{k\to \infty} q^{\frac{(\rho_k,\rho_k)}{2}}
\Bigl(\ga_x^{-1}{\mathfrak{P}}\Bigr)
(q^{x-\rho_k})=\\
&\lim_{k\to \infty} q^{\frac{(\rho_k,\rho_k)}{2}}
q^{-\frac{(x-\rho_k,x-\rho_k)}{2}}\,
{\mathfrak{P}}
(q^{x-\rho_k})=\lim_{k\to \infty}\,
q^{(x\,,\,\rho_k)}\,
\mathfrak{P}(q^{x-\rho_k}).\notag
\end{align}

In this definition, $\La$ remains untouched, so the limit 
is a $W$\~invariant function with respect to $\La$.
As a matter of fact, the key fact we need is
the existence of the limit
\begin{align}\label{whitpsi}
& \lim_{k\to \infty} 
\Psi(q^{x-\rho_k},\La;q,t)\ =\ \widetilde{\Psi}(X,\La;q)
\end{align}
for $\Psi(X,\La;q,t)$ from (\ref{pexla}). Let us calculate
the Whittaker $\widetilde{\Psi}$ in full detail. It is
essentially {\em a generating function for the 
$\overline{P}$\~polynomials}\,; see Proposition
\ref{BAREP} and related formulas for the definition of these
polynomials. 

\begin{theorem}\label{WHITG}
(i) Provided that $|q|<1$,
the Whittaker function $\widetilde{\mathfrak{P}}^\circ_x$ is given
by the formula 
\begin{align}
&\widetilde{\mathfrak{P}}^\circ_x(X,\La)
\tga_x^{\ominus}\tga_\la^{\ominus}\notag\\
=\ \widetilde{\Psi}(X,\La;q)&
\equal \sum_{b\in B_-}\, q^{\frac{(b,b)}{2}} 
\frac{X_{b_+}\ \overline{P}_b(\La^{-1}) } {
\prod_{i=1}^n\prod_{j=1}^{-(\al_i^\vee,b)}
(1-q_i^{j})
}\,,
\label{pexlabar}
\end{align}
where the power series in the right-hand side
is well defined coefficient-wise and 
converges everywhere; see
(\ref{limpct}) and (\ref{normppolsbar}). The formula
for $\widetilde{\mathfrak{P}}_x$ is with 
$\ga_x^{-1}\ga_\la^{-1}$ instead of 
$\tga_x^{\ominus}\tga_\la^{\ominus}$ and with the same
summation in the right-hand side.

(ii) The ratio of the functions 
$\widetilde{\mathfrak{P}}^\circ_x(X,\La)$,
$\widetilde{\mathfrak{P}}_x(X,\La)$ is $B$\~periodic 
with respect to $X$ and $\La$. The dependence on $\La$
is governed by (\ref{LfLa}) for the limits 
$\overline{L}_{a_+}^\La$  of the operators 
$L_{a_+}$ as $t\to 0$ upon $X\mapsto \La$\,:
\begin{align}
&\overline{L}_{a_{{}_+}}^\La\,(\widetilde{\mathfrak{P}}_x(X,\La))
\ =\ X_{a_{{}_-}}^{-1}\, \widetilde{\mathfrak{P}}_x(X,\La),\ 
X_{a_{{}_-}}^{-1}=X_{-w_{{}_0}(a_{{}_+})}.
\label{LfLadual}
\end{align}
In terms of $X$, these
functions satisfy the $q$\~Toda system of difference equations:
\begin{align}
&\widetilde{L}_{a_+}
(\widetilde{\mathfrak{P}}_x(X,\La))\ =\ 
(\,\sum_{a'\in W(a_+)}\,\La_{a'}^{-1})\  
\widetilde{\mathfrak{P}}_x(X,\La),
\label{LfLaW}\\
&\widetilde{L}_{a_+}\equal
\lim_{t\to 0}\,q^{-(a_+\,,\,\rho_k)}
\Bigl(q^{(x\,,\,\rho_k)}(\Ga_{\rho_k}^{-1}\,L_{a_+}\,\Ga_{\rho_k})
q^{-(x\,,\,\rho_k)}\Bigr),\label{qTodaop}\\ 
&\Ga_{b}(F(X))=F(q^{b}X), 
\ \Ga_b X_a=q^{(b,a)}X_a\Ga_b \for b\in \C^n. 
\notag
\end{align}
Here the difference operators $L_{a_+} (a_+\in B_+)$
from (\ref{Lf}) are conjugated by the translation  
$\Ga_{-\rho_k}$ (it is $(\rho_k)'$ in the notation
from (\ref{ondon})) and then by the operator of multiplication by
$q^{(x\,,\,\rho_k)}$. \sq
\end{theorem}

We note that $X_{b_+}\, \overline{P}_b(\La^{-1})$
in the summation for $\widetilde{\Psi}$
can be replaced by
$X_{b}^{-1}\, \overline{P}_b(\La)$.
Recall that, generally, $P_b(X^{-1})=P_{\varsigma(b)}(X)$
and $b\mapsto \varsigma(b)=-w_0(b)$ does not change
the coefficients in the summation from 
(\ref{pexlabar}). The formulas for operators
$\widetilde{L}_{a_+}$ are simple to calculate for
{\em minuscule} $a_+$; see \cite{Et1} ($A_n$) and
the rank one case below for examples.

Concerning the notation, one can introduce the
functions  $\widetilde{\mathfrak{P}}^\circ_\la$,
$\widetilde{\mathfrak{P}}_\la$
using the Whittaker limits 
with $\La,\la$ instead of $X,x$, but we do not need
these functions in the paper. Nevertheless, we put $x$ 
in $\widetilde{\mathfrak{P}}_x$ (not always)
to emphasize that the Whittaker
limit makes the dependence on $X$ and $\La$ {\em asymmetric}. 
\smallskip

The construction of the Toda operators in terms of
the Macdonald operators (and their various degenerations)
is essentially due to Inozemtsev and Etingof. 
The paper \cite{Et1} contains a systematic consideration
of the difference case.. This
paper is mainly about $GL_n$, but our (\ref{qTodaop})
is quite analogous to the limiting procedure there,
as was expected in Remark 1 at the end of \cite{Et1}.

We remark that our $q$\~Toda operators are ``dual" to those from 
\cite{Et1,GLO1} (the translation operators must be replaced 
by their inverses), which is connected with
our choice of the limit $t\to 0$ versus $t\to\infty$ in
these papers. The relation will be discussed below in 
greater detail.
\smallskip

\begin{theorem}\label{SHINT}
Continuing the previous theorem, let $X=q^c$ for $c\in B_-\,$.
Then the Shintani-type identity holds:
\begin{align}\label{shintq}
&\tga^{\ominus}(1)\,\widetilde{\mathfrak{P}}_x^\circ(q^c,\La)\ =\
\overline{P}_c(\La) 
\prod_{i=1}^n\prod_{j=1}^{\infty}\bigl(
\frac{1}{1-q_i^j }\bigr),
\end{align}
where $\tga^{\ominus}(1)=\sum_{b\in B} q^{b^2/2}$.
More explicitly, 
\begin{align}\label{shintqg}
&\sum_{b\in B_-} 
\frac{q^{(c-b,c-b)/2}\,  \overline{P}_b(\La) } {
\prod_{i=1}^n\prod_{j=1}^{(\al_i^\vee,\,b_+)}\,
(1-q_i^{j}) }\notag\\
&=\ \tga^{\ominus}(\La)\overline{P}_c(\La)
\prod_{i=1}^n\prod_{j=1}^{\infty}\bigl(
\frac{1}{1-q_i^j }\bigr). 
\end{align}
\end{theorem}
{\em Proof.} Due to (\ref{whitgen}),
\begin{align}\label{whitpr1}
&\tga^{\ominus}(q^c)\widetilde{\mathfrak{P}}^\circ_x(q^c,\La)
=\ \lim_{t\to 0} \frac{
\tga^{\ominus}(q^{c-\rho_k}) }{
\tga^{\ominus}(q^{ -\rho_k})}\,
\Bigl(\mathfrak{P}^\circ (q^{c-\rho_k},\La)\Bigr)\,.
\end{align}

Applying the identity (\ref{hatmuxsym}) for
$X$ transposed with $\La$ (the duality) inside
$\Bigr(\cdot\Bigl)$,
\begin{align}\label{whitpr2}
&\tga^{\ominus}(q^c)\widetilde{\mathfrak{P}}^\circ_x(q^c,\La)  
=\lim_{t\to 0} \frac{\tga^{\ominus}(q^{c-\rho_k})}
{\tga^{\ominus}(q^{ -\rho_k})}
\Bigl(\frac{P_c(\La)}{P_c(q^{-\rho_k})} 
\prod_{i=1}^n\prod_{j=1}^{\infty}
\frac{1}{1-q_i^j}
\Bigr)\,.
\end{align}

Recall that
\begin{align}
&\langle \ga\rangle_{\rho_k}\ =\ 
\langle \ga\rangle_{c-\rho_k}\ =\ 
\tga^{\ominus}(q^{c-\rho_k})q^{\frac{(c-\rho_k,c-\rho_k)}{2}},
\label{lgarhoky} 
\end{align}
where we use that $c$ is from $B$; see (\ref{lgarhok}). Hence,
\begin{align}\label{periwit}
&\frac{\tga^{\ominus}(q^{c-\rho_k})}
{\tga^{\ominus}(q^{ -\rho_k})}=q^{(c,\,\rho_k)-c^2/2},
\end{align}

Moving $q^{(c,\rho_k)}$ from (\ref{periwit})
to the denominator and
combining it with $P_c(q^{-\rho_k})$,
we apply (\ref{limpct}):
$$
\lim_{t\to 0}q^{-(c,\rho_k)}P_c(q^{-\rho_k})\ =\ 1.
$$

Finally, we move  $q^{-c^2/2}$ from (\ref{periwit}) to
the left-hand side of (\ref{whitpr2}) and observe that 
$q^{c^2/2}\tga^{\ominus}(q^c)$ does not depend on $c$,
so it equals $\tga^{\ominus}(1)$. \sq
\smallskip

We note that by making $q=0$ in (\ref{shintqg}), 
we arrive at the trivial
identity $\overline{P}_b(\La;q=0)=\overline{P}_b(\La;q=0)$,
where $\overline{P}_b(\La;q=0)$ is the classical character
for the dominant weight $w_0(b)$.

The $p$\~adic limit $q\to 0$ (in this setting)
transforms (\ref{hatmuxsym}) to the classical
Shintani-Casselman-Shalika formulas. 
See \cite{C101} concerning the $p$\~adic degeneration
of the DAHA theory (the limit $q\to\infty$ is considered
there).  
\medskip

\subsection{One-dimensional theory}
We will begin with the explicit formula for
the $\overline{P}$\~polynomials in the case
of $A_1$. The formulas for the Rogers polynomials
are well known as well as for their limits as $t\to 0$.
Such limits are the {\em $q$\~Hermite polynomials}
introduced by Szeg\"o and considered in many works; see, e.g.,
\cite{ASI}.
Let us re-establish the formulas we need
for these polynomials using the (nonsymmetric)
intertwining operators.

Let $\al=\al_i=\vth$, $s=s_1$, $\om=\om_1=\rho$; so
$\al=2\om$ and the standard invariant form is
$(n\om,m\om)=nm/2$.
Similarly, 
$$
X=X_\om=q^x,\ X(q^{n\om})=q^{n/2},\ 
\Ga(F(X))=F(q^{1/2} X), 
$$
i.e., $x(n\om)=n/2$, $\Ga(x)=x+1/2$, $\Ga X=q^{1/2} X\Ga$.

We will also use $\pi\equal s\Ga: X\mapsto q^{1/2}X^{-1}$;
then $\pi^2=$id\, and $Y=Y_{\om}=\pi T$ in DAHA of type $A_1$.
Concerning the Gaussians, $(x,x)=x_{\al}x_{\om}=2x^2$ and
$\ga=q^{(x,x)/2}=q^{x^2}$; note that
$\ga(q^{n\om})=q^{n^2(\om,\,\om)/2}=q^{n^2/4}$. Also,
$x-\rho_k=\Ga^{-k}(x)=x-k/2$ and $q^{(x,\rho_k)}=q^{xk}$
in the formulas for the Whittaker limit.
Here $t=t_0=q^k \for k\in \C$.

In the limit $t\to 0$, $\overline{T}=T(t=0)$ and  
$$
\overline{P}_n=\overline{P}_{-n\om}=
\overline{E}_{-n\om}.
$$
For instance, $\overline{P}_0=1$,\ $\overline{P}_1=X+X^{-1},$
\begin{align}
\overline{P}_2=&X^2+X^{-2}+1+q,\ 
\overline{P}_3=X^3+X^{-3}+\frac{1-q^3}{1-q}(X+X^{-1}),\notag\\
&\overline{P}_4=X^4+X^{-4}+\frac{1-q^4}{1-q}(X^2+X^{-2})
+\frac{(1-q^4)(1-q^3)}{(1-q)(1-q^2)}.\label{overpsmall}
\end{align}
Generally, for the monomial symmetric functions
$M_0=1$, $M_n=X^n+X^{-n}$ for $n>1$, 
\begin{align}
&\overline{P}_n=M_n+\sum_{j=1}^{[n/2]}\,
\frac{(1-q^{n})\cdots(1-q^{n-j+1})}
{(1-q)\cdots(1-q^{j})}\, M_{n-2j}.\label{overpgen}
\end{align}

The norm formulas from
(\ref{mubar}), (\ref{constermbar}), (\ref{normppolsbar}) 
read as follows:
\begin{align}\label{normppolsbara1}
&\lan \overline{P}_{m}(X)\overline{P}_{n}(X)
\overline{\mu}_\circ\ran\ 
=\ \de_{mn}\prod_{j=1}^{n}
(1-q^{j})\,,
\end{align}
where $m,n=0,1,\ldots$,\  $\overline{\mu}_\circ=
\overline{\mu}/\lan\overline{\mu}_\circ\ran$ for
the classical theta-function
\begin{align}
&\overline{\mu} =  
\prod_{j=0}^\infty (1-X^2 q^{j})
(1-X^{-2}q^{j+1}),\ 
\langle\overline{\mu}\rangle = 
\prod_{j=1}^{\infty} \frac{1}
{1-q^j}.
\label{constermbar1}
\end{align} 

Due to Theorem \ref{PHIEB}, the composition
$\r=(1+\overline{T})X\pi$ is the raising operator
for the $\overline{P}$\~polynomials. Namely,
upon the restriction, Red, to the symmetric polynomials:  
\begin{align}\label{raisingbar1}
&q^{\frac{n}{2}}R(\overline{P}_n)=\overline{P}_{n+1},
\where R=\hbox{Red}(\r)=\frac{X^2\Ga^{-1}-X^{-2}\Ga}
{X-X^{-1}}.
\end{align}
This readily gives (\ref{overpgen}).
\smallskip

{\bf Rogers polynomials.} The  
counterparts of these formulas for the Rogers polynomials
are well-known
(see, e.g., \cite{AI} and \cite{C101}, Chapter 2). Let us 
list them for the sake of completeness.
First,
\begin{align}\label{constermbar1t}
&\mu\ =\   
\prod_{j=0}^\infty \frac{(1-X^2 q^{j})
(1-X^{-2}q^{j+1})}{(1-X^2 tq^{j})
(1-X^{-2}tq^{j+1})},\where  \\
&\langle\mu\rangle\ =\ \hbox{Constant Term\,}(\mu)\ =\  
\prod_{j=1}^{\infty} \frac{(1-tq^j)^2}
{(1-t^2q^j)(1-q^j)}.\notag
\end{align} 
Switching to $\mu_\circ=\mu/\lan \mu\ran,$  
\begin{align}\label{normppolsbara1t}
&\lan P_{m}(X)P_{n}(X)
\mu_\circ\ran\ 
=\ \de_{mn}\prod_{j=0}^{n-1}\frac{(1-q^{j+1})
(1-t^2q^{j})}{(1-tq^{j+1})
(1-tq^{j})}\,,
\end{align}
as $m,n=0,1,2,\ldots\,.$ The explicit formulas are as follows:
\begin{align}
&P_n=M_n+\sum_{j=1}^{[n/2]}\,
\prod_{i=0}^{j-1}
\frac{(1-q^{n-i})}
{(1-q^{1+i})}\,
\frac{(1-tq^{i})}{(1-tq^{n-i-1})}
\, M_{n-2j}.\label{overpgent}
\end{align}
\medskip

{\bf The $L$\~operators.} We will begin with the
formula for the $q$\~Toda operator
from (\ref{qTodaop}):
\begin{align}\label{qtoda1}
\widetilde{L}=
&\lim_{t\to 0}\, \bigl(q^{kx}\,\Ga_{k}^{-1}\,L\,\Ga_{k}\,
q^{-kx}\,\bigr)=(1-X^{-2})\Ga+\Ga^{-1}
\end{align}
for $\Ga_k X\equal t^{k/2}X\Ga_k\,$ and
the well-known operator 
\begin{align}\label{Macdop}
L=L_{\om}=\hbox{Red}(Y+tY^{-1})=
\frac{1-tX^{2}}{1-X^2}\Ga+\frac{1-tX^{-2}}{1-X^{-2}}\Ga^{-1},
\end{align} 
diagonalizable in terms of Rogers' polynomials.
We will also use:
\begin{align}\label{qtodaga1}
&\widetilde{L}_\ga=\ga^{-1}\widetilde{L}\ga=
q^{1/4}\bigl(X\Ga+X^{-1}(\Ga^{-1}-\Ga)\bigr)
\for \ga=q^{x^2}.
\end{align}

For the straightforward specialization of $L$
at $t=0$, one has: 
\begin{align}
&\overline{L}\ =\ \lim_{t\to 0}L\ =\ 
(1-X^{2})^{-1}\Ga+(1-X^{-2})^{-1}\Ga^{-1},\notag\\ 
&\overline{L}_\ga=\ga^{-1}\overline{L}\ga=
-q^{1/4}(X-X^{-1})^{-1}(\Ga-\Ga^{-1}).
\label{qtlga1}
\end{align}
The latter operator is proportional to  
the so-called Askey -Wilson divided difference operator,
which serves as the {\em shift operator}
in the theory of Rogers' polynomials (with any $t$)
and the basic hypergeometric function. See 
\cite{AI} and also \cite{C101}, Chapter 2. 

Its defining property is the relation
\begin{align}\label{shifta1}
&\overline{L}_\ga(\overline{P}_n)\ =\ 
-q^{1/4}(q^n-q^{-n})\overline{P}_{n-1},\ 
n=1,2\ldots\,.
\end{align}
Let us give a convenient reference concerning
(\ref{raisingbar1}),(\ref{shifta1}):
\cite{OS}, formulas (20-25). 
\medskip

\subsection {Whittaker function for
\texorpdfstring{{\mathversion{bold}$A_1$}}{A1}}
Provided that $|q|<1$, we can now introduce the Whittaker 
function $\widetilde{\mathfrak{P}}^\circ_x$
from the relation: 
\begin{align}
&\widetilde{\mathfrak{P}}^\circ_x(X,\La)
\tga_x^{\ominus}\tga_\la^{\ominus}\ =\ \widetilde{\Psi}(X,\La)
\equal \sum_{n=0}^\infty 
q^{\frac{n^2}{4}} 
\frac{X^n\ \overline{P}_n(\La) } {
\prod_{j=1}^{n}
(1-q^{j})
}\,,
\label{pexlabara1}
\end{align}
where $\tga_x^{\ominus}=\sum_{j=-\infty}^\infty q^{j^2/4}X^j$
($\tga_\la^{\ominus}$ is defined in terms of $\la$).

The function $\widetilde{\Psi}(X,\La)$ is actually 
the generating function for $q$\~Hermite polynomials.
It is directly connected with the, so-called, quadratic 
$q$\~exponential function;
see \cite{Sus} formulas (26),(27) and the
references there. 
Its interpretation as a $q$\~Whittaker function
(upon the multiplication by the Gaussians)
does not seem to have been noticed, although 
the difference equation for $\widetilde{\Psi}(X,\La)$ 
was certainly known (formula (19) ibid.). The one-dimensional
Shintani-type formulas, (\ref{shintqa1}) below and especially
its $q,t$\~generalization,
seem new. Some related formulas like (\ref{pbgauss}) can
be deduced from known identities (at level of ${}_6\Psi_6$);
the multidimensional theory is new.

The power series $\widetilde{\Psi}(X,\La)$ converges everywhere.
The $\La$\~dependence (see (\ref{LfLadual})) readily
follows from (\ref{shifta1}):
\begin{align}
&\overline{L}_\ga^\La\ (\widetilde{\Psi}(X,\La))
\ =\ X\, \widetilde{\Psi}(X,\La),\ \overline{L}_\ga^\La=
\overline{L}_\ga(X\mapsto \La).
\label{LfLaduala1}
\end{align}

In terms of $X$, the function $\widetilde{\Psi}$ satisfies the
{\em $\ga$\~twisted} $q$\~Toda equation, which 
reads as follows:
\begin{align}
&\widetilde{L}_\ga\,
(\widetilde{\Psi}(X,\La))\ =\ 
(\La+\La^{-1})\,\widetilde{\Psi}(X,\La).
\label{LfLaWa1}
\end{align}
\medskip

{\bf The $q,t$\~case.} The formula for the 
global spherical function $\mathfrak{P}^\circ$ reads as follows:
\begin{align}\label{pexlabara1t}
&\mathfrak{P}^\circ(X,\La)
\,\tga_x^{\ominus}\tga_\la^{\ominus}/\tga^{\ominus}(q^k)\ =\ 
\Psi(X,\La)\\
&\equal \sum_{n=0}^\infty 
t^{\frac{n}{2}}q^{\frac{n^2}{4}} 
P_n(X) P_n(\La) 
\prod_{i=0}^{n-1}\frac
{(1-tq^{i})(1-tq^{i+1})}{(1-t^2q^{i})(1-q^{i+1})}\,.
\notag
\end{align}
We come to a variant of the basic hypergeometric function.
\medskip

{\bf Shintani-type formula}.
Let us consider Theorem \ref{SHINT} in the $A_1$\~case;
we plug in $X=q^{-n/2}$ for $n=0,1,\ldots,$\,. Then 
\begin{align}\label{shintqa1}
&q^{n^2/4}\,\widetilde{\Psi}(q^{-n/2},\La)\ =\
\tga^{\ominus}(\La)\overline{P}_n(\La) 
\prod_{j=1}^{\infty}\Bigl(
\frac{1}{1-q^j }\Bigr).
\end{align}
Recall that $\tga^{\ominus}(q^\la)=\sum_{j\in \Z} q^{j\la+j^2/4}$.
Here the left-hand side and the right-hand side
coincide as Laurent series or as analytic functions.

This formula becomes a trivial identity for $q=0$, i.e.,
in the case of the classical characters 
$$
\overline{P}_n(X;q=0)\ =\ 
\frac{X^{n+1}-X^{-n-1}}{X-X^{-1}}.$$
\medskip

\subsection{The case  
\texorpdfstring{{\mathversion{bold} $|q|>1$}}{of large |q|}}

Generally, the Whittaker-type
limiting procedure as $t\to \infty$
is naturally connected with the theory at $|q|>1$
and can lead to new formulas.
However, in the symmetric setting of this paper, 
there is a direct connection between the Whittaker functions
defined for $|q|<1,t\to 0$ and $|q|>1,t\to \infty$, 
which we are going to discuss now. 
\smallskip

We follow \cite{C5}
and use $\tga^{\oplus}$ instead of $\tga^{\ominus}$
and $\ga$ instead of $\ga^{-1}$. In the 
nonsymmetric setting, the corresponding
global spherical function is really different from
that for $|q|<1$. However, there exists a simple
connection in the symmetric case.

The $q,t$\~definition we need is as follows (cf. (\ref{pexla})):
\begin{align}
& 
\tga_x^{\oplus}\tga_\la^{\oplus}\,
\mathfrak{P}_\star^\circ/\tga^{\oplus}(q^{\rho_k})\  =\ 
\sum_{b\in B_-} q^{-\frac{(b,b)}{2} +(\rho_k,b)} 
 \frac{P_b(X)\ P_b(\La)} {
\langle P_b(X) P_b(X^{-1})\mu_\circ\rangle},
\label{pexla1}
\end{align}
where $q>1$ and $\mathfrak{P}_\star^\circ$ satisfies
the claims of Theorem \ref{GLOBSPH}.
The Whittaker limiting procedure requires here
taking $t\to \infty$  for ensuring the convergence. 

The formulas are:
\begin{align}\label{whitgen1}
& \tga_\la^{\oplus}\tga_x^{\oplus}
\widetilde{\mathfrak{P}}_\star^\circ 
=\tga_\la^{\oplus}
\lim_{t\to \infty} \frac{\tga_x^{\oplus}(q^{x-\rho_k})}
{\tga_x^{\oplus}(q^{\rho_k})}\,
\mathfrak{P}^{\circ}_\star
(q^{-\rho_k}X)\\
&\equal\widetilde{\Psi}_\star(X,\La;q)\ =\ 
\sum_{b\in B_-} q^{-\frac{(b,b)}{2}} 
\frac{X_{b}\ P_b(\La;q,t\to \infty) } {
\prod_{i=1}^n\prod_{j=1}^{-(\al_i^\vee,b)}
(1-q_i^{-j})}.
\label{pexlabar1}
\end{align}
Cf. (\ref{whitgen}) and (\ref{pexlabar}). Here 
\begin{align}\label{ptinvlim}
&\lim_{t\to\infty}P_b(\La;q,t)=
\lim _{t^{-1}\to 0}P_b(\La^{-1};q^{-1},t^{-1})=
\overline{P}_b(\La^{-1};q^{-1}).
\end{align}
Therefore $\widetilde{\Psi}_\star(X,\La; q)$ 
simply {\em coincides}
with $\widetilde{\Psi}(X^{-1},\La; q^{-1})$ in the notation
from (\ref{pexlabar}). 

We conclude that $\widetilde{\mathfrak{P}}_\star^\circ$
satisfies the eigenvalue problem
\begin{align}
&\widetilde{L}^{\star}_{a_+}
(\widetilde{\mathfrak{P}}_\star(X,\La))\ =\ 
(\,\sum_{a'\in W(a_+)}\,\La_{a'}^{-1})\  
\widetilde{\mathfrak{P}}_\star(X,\La),
\label{LfLaW1}\\
&\widetilde{L}^{\star}_{a_+}\equal
\lim_{t\to \infty}\,q^{-(a_+\,,\,\rho_k)}
\Bigl(q^{-(x\,,\,\rho_k)}(\Ga_{\rho_k}^{-1}\,L_{a_+}\,\Ga_{\rho_k})
q^{(x\,,\,\rho_k)}\Bigr),\label{qTodaop1}
\end{align}
where $\Ga_{b}(F(X))=F(q^{b}X)$. Compare with  
(\ref{qTodaop}); the conjugation by $q^{(x\,,\,\rho_k)}$
there is replaced by the conjugation by $q^{-(x\,,\,\rho_k)}$.
Thus the operators $L_{a_+}^\star(X,q^{-1})$ generalize
those considered in \cite{Et1,GLO2}. For instance,
in the one-dimensional case 
in the notation from (\ref{qtoda1}):
\begin{align}\label{qtodaa1}
\widetilde{L^\star}&=
\lim_{t\to \infty}\, \bigl(q^{-kx}\,\Ga_{k}^{-1}\,L\,\Ga_{k}\,
q^{kx}\,\bigr)\\
=\lim_{t\to \infty}
t^{-1/2}&\Bigl(
\frac{t(t^{-1}X^2)-1}{t^{-1}X^2-1}(t^{1/2}\Ga)+
\frac{t(tX^{-2})-1}{tX^{-2}-1}(t^{-1/2}\Ga^{-1})\Bigr)\notag\\
&=(1-X^{2})\Ga+\Ga^{-1}.\notag
\end{align}
\medskip

\setcounter{equation}{0}
\section{Harmonic analysis topics}
The real integration or Jackson integration is, generally,
necessary when the Gaussian $\ga^{-1}$ in the constructions
above is replaced by $\ga$. A typical example is
as follows. Let us consider the DAHA-Fourier transform 
in terms of the constant term functional
(or using the imaginary integration) in the space of Laurent
polynomials multiplied by $\ga^{-1}$. Then 
the inverse transform will involve the Jackson (or real)
integration and the proper choice of the
Gaussian is $\ga$ instead of $\ga^{-1}$. 

Such ``switch" of the Gaussians
is necessary algebraically due to the properties of
the involution of DAHA that governs the Fourier transform.
Correspondingly, the contour of integration, real
or imaginary,  must ensure the convergence, i.e., its
choice is of analytic nature. The direction is real 
for $\ga$ and imaginary for $\ga^{-1}$.
It is of course for $|q|<1$; if $|q|>1$ then it must be 
the other way round. Generally,
especially, in the absence of the Gaussians 
(for instance, in the Harish-Chandra theory), the
directions, real or imaginary, are selected to match the 
growth estimates for the spherical function, 
used as kernels of the corresponding transforms. 

We establish such estimates in the {\em real direction}. 
The theory appears surprisingly 
``precise", although the results of the paper are far from
being complete. Only the first term of
the asymptotic expansion is obtained. We note that in 
our setting, the global
spherical function is periodic in the imaginary direction, 
so the imaginary growth estimates are irrelevant.
We stick to the Jackson
integration, which is actually very similar to 
the ``classical" case of real integration; the estimates we
obtain serve both theories.
\smallskip

\subsection{Growth estimates}
It is possible to evaluate the growth of the
global $q,t$\~spherical function $\mathfrak{P}(X,\La;q,t)$
from Theorem \ref{GLOBSPH} in the real directions. 
Let $0<q<1$, $t_\nu=q_\nu^{k_\nu}$ (or, simply, $t=q^k$)
for $k_\nu\in \C$ provided the existence of all 
spherical symmetric polynomials $\{\p_{b_-}\}$, 
equivalently, provided that
the polynomial representation is semisimple and the radical
of the evaluation pairing vanishes (see \cite{C103}).
The assumption $\Re k_\nu>-1/h_\nu$ for the 
{\em Coxeter numbers}
$h_\nu=1+(\rho,(\th')^\vee)$, where $\th'=\th,\vth$
for $\nu=\nu_{\lng},\nu_{\sht}$, is sufficient (but
not necessary).

For $x\in \C^n$, let $x_+\equal u(x)$ where $u(\Re(x))$ 
is a unique vector belonging to the closure 
$\overline{\mathfrak{C}}_+= \sum_{i=1}^n\, \R_+\,\om_i$ 
of the standard positive
nonaffine Weyl chamber $\mathfrak{C}_+=
\sum_{i=1}^n \R_{>0}\,\om_i$.

Given a $p$\~sequence of vectors 
$\mathbf{x'}=\{x'_1,\ldots, x'_p\}\subset \R^n$
and a $p$\~sequence of positive integers 
$\mathbf{n}=\{n_1,\ldots, n_p\}$, 
we use the dot-notation 
$\mathbf{n}\cdot \mathbf{x'}$ for
$\sum_{j=1}^p n_j x'_j$.

\newcommand{\bla}{\mbox{\mathversion{bold}$\la'$}}
\newcommand{\nx}{\mathbf{n}\cdot \mathbf{x'}}
\newcommand{\ny}{\mathbf{n}\cdot \mathbf{y'}}
\newcommand{\nbla}{\mathbf{n}\cdot \bla}

\begin{theorem}\label{GROWP}
(i) For arbitrary $x,\la\in \C^n, k_\nu\in \C$, we set\,: 
\begin{align}\label{pdagcirc}
&\mathfrak{P}^\circ_{\dagx}(x,\la;q,k)\equal 
\frac{\tga^\ominus(q^x)\tga^\ominus(q^{\la})}
{\tga^\ominus(q^{x_+ +\la_+ -\rho_k})\tga^\ominus(q^{\rho_k})}
\,\mathfrak{P}^\circ(q^{x},q^\la; q,q^k).
\end{align}
Given a real $p$\~sequence
$\mathbf{x'}$, let the components of 
$\mathbf{n}$ tend to $+\infty$ in an arbitrary
way provided that 
$
(\mathbf{n}\cdot \mathbf{x'})_+\in 
\mathfrak{C}_+.
$
Then the limit
\begin{align}\label{growthe}
&\lim_{\mathbf{n}\to\infty}
\mathfrak{P}^\circ_{\dagx}
(x+\mathbf{n}\cdot \mathbf{x'},
\la;q,k)
\end{align}
exists if $\Re(\la)_+\in \mathfrak{C}_+$; moreover,
it depends only on $\la$ and is nonzero for all such $\la$.
Here we choose $x$ to ensure that 
$\tga^\ominus(q^{(x+\mathbf{n}\cdot \mathbf{x'})_+
+(\la)_+ -\rho_k})
\neq 0$ for any $\mathbf{n}$.

Under the same constraints, consider
$\mathfrak{P}^\circ_{\dagx}
(x+\mathbf{n}\cdot \mathbf{x'},
\la+\mathbf{n}\cdot \hbox{\mathversion{bold}$\la'$};q,k)$
for a real $p$\~sequence {\mathversion{bold}$\la'$} satisfying
$
(\mathbf{n}\cdot \hbox{\mathversion{bold}$\la'$})_+\in 
\mathfrak{C}_+\,.
$
Then the limit exists too and is an absolute nonzero
constant depending only on $q,k$. 

(ii) In the case of the
Whittaker function $\widetilde{\mathfrak{P}}$, 
we remove $k$ from the formulas and 
replace $x_+$ by $-x$: 
\begin{align}\label{wtildedag}
&\widetilde{\mathfrak{P}}^\circ_{\dagx}(x,\la; q)\equal
\frac{\tga^\ominus(q^x)\tga^\ominus(q^{\la})}
{\tga^\ominus(q^{\la_+ -x})} 
\,\widetilde{\mathfrak{P}}^\circ(q^{x},q^\la; q)\,.
\end{align}
Provided that $\mathbf{n}\cdot\mathbf{x'}\in -\mathfrak{C}_+$
(it was not needed in the $q,t$\~case), the 
claims from (i) hold true for  
\begin{align}\label{growthew}
&\lim_{\mathbf{n}\to\infty}\widetilde{\mathfrak{P}}^\circ_{\dagx}
(x+\mathbf{n}\cdot \mathbf{x'},
\la+\mathbf{n}\cdot \hbox{\mathversion{bold}$\la'$};q).
\end{align}
Here $\tga^\ominus(q^{\{\cdot\}})$ is nonzero at
$(\la+\mathbf{n}\cdot 
\hbox{\mathversion{bold}$\la'$})_+
-(x+\mathbf{n}\cdot \mathbf{x'})\,;$
we continue to assume that 
$
(\mathbf{n}\cdot \mathbf{x'})_+\in 
\mathfrak{C}_+
$
and, correspondingly, either  $\Re(\la)_+\in \mathfrak{C}_+$ for
{\mathversion{bold}$\la'$}$=0$ or
$
(\mathbf{n}\cdot \hbox{\mathversion{bold}$\la'$})_+\in 
\mathfrak{C}_+\,.
$
\sq
\end{theorem}

Let us comment on the proof. 
The justification of (i) involves the analysis
of the corresponding difference equations for
$\mathfrak{P}$ in the limit of large $x$ and/or 
large $\la$, but we use the explicit formulas too.
We note that the {\em asymptotic} difference equations provide
the asymptotic limit (the factor in the definition of 
$\mathfrak{P}^\circ_{\dagx}$ from (\ref{wtildedag})) only up to 
a periodic function. So we need to use that both, 
$\mathfrak{P}^\circ$
and $\mathfrak{P}^\circ_{\dagx}$, are meromorphic.

Using the asymptotic differential or difference equations for 
the analysis of the limiting behavior of solutions
of the corresponding equations is standard. 
Given $x,\la$, we restrict ourselves with the bi-lattice  
$\{x+P_+,\la+P_+\}$
and evaluate the values of $\mathfrak{P}^\circ$
there step-by-step using the corresponding difference
equations. 

For instance, we treat the $L$\~operator from
(\ref{Macdop}) in the case of $A_1$ as a recurrence
for calculating the value at $x+(n+1)\om_1$ in terms of
the values at $x+n\om_1$ and $x+(n-1)\om_1$,
where the coefficients tend to constants
as $n\to\infty$. The stabilization is of exponential 
type, which simplifies the necessary estimates.
In the $\la$\~space, we use the
$x\leftrightarrow\la$\~duality and the $\la$\~counterpart
of the $L$\~operator. We note that in the rank one case, 
there is the classical theory by Birkhoff, which was developed
recently in several works. 

The multi-dimensional case
is not very different, as far as Theorem \ref{GROWP}
is concerned. Similar problems were considered in the theory
of Quantum Knizhnik-Zamolodchikov equations. By the way,
the equivalence of the QAKZ and the eigenvalue problem
under consideration (see \cite{C101}), generally, can be
used here. Our approach is direct.  For arbitrary root
systems, the formulas for the $L$\~operators are not
explicit, but we need only the stabilization estimates for
their coefficients. This approach becomes significantly more
transparent in the nonsymmetric theory, where we can use
directly the intertwining operators instead of the 
difference relations. The nonsymmetric {\em global}
functions have various symmetries including the 
transformation formulas under the action of the intertwining
operators in both, the $x$\~space and the $\la$\~space. 


Part (ii) is obtained as a limit of (i).
Taking $\mathbf{x'}$ and $\bla$ {\em real vectors} is,
actually,  
insignificant in the theorem. Since  $\mathfrak{P}^\circ$ and
$\mathfrak{P}^\circ_{\dagx}$ are $2\pi i \log(q) P^\vee$\~periodic
in the imaginary direction, it suffices to impose
the conditions from (i,ii) for their real parts only. 

We note, that due to the claim that the limits do not depend on the 
particular way the integers $\{n_i\}$ approach the infinity,
one can try to use the Shintani-type formulas. It is
assuming that the uniqueness of
$\mathfrak{P}^\circ$ is known in the corresponding analytic class 
of solutions of the spherical eigenvalue problem satisfying
the Shintani-type formulas, i.e., among the solutions
that ``go through" the $\p$\~polynomials.
To avoid misunderstanding, let
us emphasize that we do not claim or use such uniqueness
in this paper. An approach to its justification we have
in mind employs the symmetry $x\leftrightarrow \la$ or the
passage to the nonsymmetric theory.


\subsection{Exact asymptotic formulas}
Let us obtain them in the (most important) case
$\bla=0$. As a matter of fact, we do not need 
(more general but less exact) Theorem \ref{GROWP},
although the limiting difference equations are 
used in the proof of the main lemma. Recall that we deal 
with the function that is given by an explicit formula; 
the asymptotic difference equations are natural and convenient 
here but can be replaced by straightforward analysis of 
$\mathfrak{P}^\circ$ based on the asymptotic theory of 
Macdonald polynomials. It is what we are going to do in 
this section.  
 
The key ingredient
is the {\em inverse} of the
positive half of the $\mu$\~function,
a direct $q,t$\~counterpart of the 
celebrated Harish-Chandra
$c$\~function \cite{HC}:
\begin{align}
&\si(X;q,t) = \prod_{\al \in R_+}
\prod_{j=0}^\infty \frac{1-t_\al X_\al q_\al^{j}
}{
1-X_\al q_\al^{j}
}.\
\label{muhalf}
\end{align}

\begin{theorem}\label{QHARISH}
(i) Provided the conditions of part (i) of
Theorem \ref{GROWP} for $\bla=0$,
including $\Re(\la)_+\in \mathfrak{C}_+$,
\begin{align}\label{growthex}
&\lim_{\mathbf{n}\to\infty}
\mathfrak{P}^\circ_{\dagx}
(x+\mathbf{n}\cdot \mathbf{x'},
\la;q,k)= \varrho(q,t)\,
\sigma(q^{\la_+};q,t)\for\\
&\varrho(q,t)\equal \lan\, \sum_{w\in W} w(\mu)\, \ran\,=\,
\lan \mu \ran\,
\prod_{\al>0}\,\frac{1- q^{(\rho_k,\al)}}
{1- t_\al q^{(\rho_k,\al)}} \notag\\
&=\ \prod_{\al>0}
\prod_{j=1}^{\infty} \frac{ (1- q^{(\rho_k,\al)+(j-1)\,\nu_\al})
(1- q^{(\rho_k,\al)+j\,\nu_\al})}
{(1-t_\al q^{(\rho_k,\al)+(j-1)\,\nu_\al})
(1-t_\al^{-1}q^{(\rho_k,\al)+j\,\nu_\al})}\,,
\label{consterminv}
\end{align}
where $\lan \mu\ran$ is the constant term of $\mu$
from (\ref{consterm}).

(ii) Correspondingly, 
imposing $\mathbf{n}\cdot\mathbf{x'}\in -\mathfrak{C}_+$
and the other conditions in the Whittaker case,
\begin{align}\label{growthexw}
\lim_{\mathbf{n}\to\infty}
\widetilde{\mathfrak{P}}^\circ_{\dagx}&
(x+\mathbf{n}\cdot \mathbf{x'},\la;q)\ =\ 
\lan \overline{\mu}\,\ran\,\sigma(q^{\la_+};q,0)\\
&=\ 
\prod_{i=1}^n \prod_{j=0}^{\infty} \,
\frac{1}
{(\,1-q_i^{j+1}\,)\,
(\,1-q_i^{\,(\la_+\,,\,\al_i^\vee)+j}\,)}\,.\notag
\end{align}
In contrast to this formula, assuming that 
$(\nbla)_+\in \mathfrak{C}_+$, the $\la$\~limit does
not depend on $x$:
\begin{align}\label{growthexwla}
&\lim_{\mathbf{n}\to\infty}
\widetilde{\mathfrak{P}}_{\dagx}
(x,\la+\nbla;q)\ =\
\lan \overline{\mu}\,\ran\ =\ 
\prod_{i=1}^n \prod_{j=1}^{\infty} \,
\frac{1}{1-q_i^{j}}\,.
\end{align}
The limit remains the same if we substitute
$x\mapsto x+\nx$ in (\ref{growthexwla})
for $\mathbf{x'}$
such that  $\Re((\nbla)_+ -\nx)\in \mathfrak{C}_+$.
\end{theorem}
{\em Proof}. It suffices to calculate
\begin{align}\label{growthexzz}
&\lim_{c_+\to\infty}\,
\mathfrak{P}^\circ_{\dagx}
(c-\rho_k,\la;q,k), \where c\in B_-\,,
\end{align}
and by $c_+\to\infty$, we mean that
$(\al_i,c_+)\to \infty$
for all $i=1,\cdots,n.$  Recall that $c_+=w_0(c)$,
where $c$ is always from $B_-$ in this calculation. 

Using the definition and formula (\ref{hatmuxsym}),
\begin{align}\label{shintqcal}
\mathfrak{P}^\circ_{\dagx}&
(c-\rho_k,\la;q,k)=
\frac{\tga^\ominus(q^{c-\rho_k})\tga^\ominus(q^{\la})}
{\tga^\ominus(q^{\la_+ - c})\tga^\ominus(q^{\rho_k})}\\
&\times 
\frac{P_c(q^{\la})}{P_c(q^{-\rho_k})} 
\prod_{\al\in R_+}\prod_{ j=1}^{\infty}\Bigl(\frac{ 
1- q_\al^{(\rho_k,\,\al^\vee)+j}}{
1-t_\al^{-1}q_\al^{(\rho_k,\,\al^\vee)+j} }\Bigr). 
\label{hatmuxsymcal}
\end{align}
The special value $P_c(q^{-\rho_k})$ is given by
(\ref{pebebs}); it is the
exponent  $q^{(\rho_k,\,c)}$ times the product term,
which will be combined  
(in the limit of large $c_+$) with the product from
(\ref{hatmuxsymcal}). The result is 
exactly $\varrho(q,t)$, the constant term of the
symmetrization of $\mu$ from \cite{M3,M2,C2}. 

We note that  $\lan \mu \ran$  was obtained in this 
calculation without any reference to its ``true" 
meaning as the constant term of $\mu$. 
It is interesting but not
very much surprising; in \cite{C101} the norm-formula
for Macdonald polynomials (including the constant
term formula) was actually deduced from 
the evaluation formula. Something similar occurs here.

Since $c\in B$ (actually $c \in B_-$), we can remove
it  from the  theta-functions 
$\tga^\ominus(q^{c-\rho_k})$ and 
$\tga^\ominus(q^{\la_+ - c})$, the multiplicators
are the same as for the Gaussians 
$q^{-(c-\rho_k)^2/2}$ and  $q^{-(\la_+-c)^2/2}$. It gives:
\begin{align}\label{shintqcal1}
&\frac{\tga^\ominus(q^{c-\rho_k})\tga^\ominus(q^{\la})}
{\tga^\ominus(q^{\la_+ - c})\tga^\ominus(q^{\rho_k})}\ =\ 
q^{(c,\,\rho_k-\la_+)}.
\end{align}
The factor $q^{(c,\,\rho_k)}$ will cancel the same term
from $P_c(q^{-\rho_k})$ (in the denominator).  
The remaining part of (i)
is taking the limit 
$$
\lim_{c_+\to\infty} q^{-(c,\,\la_+)}P_c(q^\la),
$$
which is a subject of the following lemma.

\begin{lemma} \label{LEMCINFTY}
Provided that $|q|<1$ and $\Re (x_+)\in \mathfrak{C}_+$,
$$
\lim_{c_+\to\infty} q^{-(c,\,x_+)}P_c(q^x)\ =\
\sigma(q^{x_+};q,t),
$$
where the limit is pointwise.
\end{lemma}
{\em Proof}. In the multiplicative notations,
$q^{(c,\,x_+)}=q^{(w^{-1}(c),\,x)}=X_{w(c)}$
for $w(\Re(x))\in \mathfrak{C}_+$, i.e., this
monomial is from the leading symmetric monomial 
function of the $P_c(X)$. Its
coefficient is $1$ by construction. One can assume
here that $w=1$ due to the $W$\~invariance of $P_c$.
Then $X_c^{-1}P_c$ will be a power series in terms
of $X_{\al_i}$ for $i=1,\cdots,n$. 

Calculating the corresponding difference equations
(in the limit of large $c_+$) is the most direct way
to identify its expansion with $\sigma(X)$.
It suffices to uses the 
leading terms of the $L$\~operators serving
the symmetric Macdonald polynomials calculated
in \cite{C2}, Proposition 3.4. Then we observe
that $\sigma(X)$ is a solution of this system
of equations. It gives the required since both are
power series  in terms of $X_{\al_i}$ with the 
constant term $1$. 
\sq
\smallskip

The lemma gives (\ref{growthex}). The Whittaker
variants from (ii) are its straightforward
limits; the condition $\Re(x)\in -\mathfrak{C}_+$ must
be imposed in (\ref{growthexw}) and no such conditions
are necessary in (\ref{growthexwla}).
Obtaining these two limits via 
the Shintani-type formulas (\ref{shintq}) seems
possible as well, however, it requires knowing that 
these formulas determine
$\mathfrak{P}^\circ$ uniquely in a proper class of 
functions, which we do not claim in this paper.

\sq
\smallskip

Lemma \ref{LEMCINFTY} is known for the
Askey-Wilson polynomials (see, e.g., \cite{Is}).
The Laurent expansion of the rank one $\mu$\~function 
is very explicit, so it is straightforward.
Paper \cite{FZ} contains a comprehensive
discussion of the $A_1$\~case.
In paper \cite{Ru}, the claim of the lemma was
obtained in the $A$\~case for the $L^2$\~convergence.
It was conjectured there (with some explicit estimate)
that the convergence is pointwise as well; see a 
discussion after formula (1.23).
Paper \cite{vD0} is an extension of \cite{Ru} to
the case of arbitrary reduced root systems 
(for the strong $L^2$\~convergence). See also \cite{vD} for the case 
of the Koornwinder polynomials (the root system $C^\vee C_n$).

Our {\em operator} approach (based on the asymptotic difference
operators) gives the pointwise convergence. We can,
generally, answer Ruijsenaars' question concerning the
pointwise estimates in compact sets. However, we will not
touch upon this (important) direction in this paper.  

As for Theorem \ref{QHARISH}, we think that 
its one-dimensional versions (for the basic hypergeometric
function or its variants) are likely to be known.
\smallskip

\subsection {The Harish-Chandra formula}
The corollary is an exact generalization of the
Harish-Chandra fundamental asymptotic formula
for the classical spherical functions. Indeed,
for $x$ approaching $\infty$ in the directions
$\mathbf{x'}$ (admissible in the sense of
Theorem \ref{GROWP}), asymptotically,

\begin{align}\label{growthexfin}
&\lim_{\mathbf{n}\to\infty}
\mathfrak{P}^\circ
(x+\mathbf{n}\cdot \mathbf{x'},
\la;q,k)\\
\sim\ \varrho(q,t)\,
&\frac{\tga^\ominus(q^{x_+ +\la_+ -\rho_k})
\tga^\ominus(q^{\rho_k})}
{\tga^\ominus(q^{x_+})\tga^\ominus(q^{\la_+})}
\,\sigma(q^{\la_+};q,t).\notag
\end{align}
Up to a simple $W$\~invariant and $B$\~periodic
factor  $C(x,\la)$, depending of course on $q,k$ (it is 
$\Z$\~periodic in terms of $k$),  
we can switch to $\mathfrak{P}$ here, replacing
all $\tga^\ominus(q^x)$ by $\ga^{-1}(q^x)=q^{-x^2/2}$. It
gives that in the limit of large $\Re(x)_+\in \mathfrak{C}_+$, 
\begin{align}\label{growthexfinx}
&\mathfrak{P}^\circ(x,\la;q,k)
\sim C(x,\la)\,\varrho(q,t)\,
q^{-(x_+,\, \la_+ -\rho_k)+
(\la_+,\,\rho_k)}\,\sigma(q^{\la_+};q,t).
\end{align}

\begin{corollary}\label{HELGJ}
We continue to assume that all spherical polynomial 
$\{\p_{b_-}\}$ exist; for instance, the conditions
$k_\nu\not\in -1/h_\nu-\Q_+$
for the Coxeter numbers $h_\nu$ of $R$ are sufficient.
Provided that  $\Re(\la)_+\in \mathfrak{C}_+$, 
the global spherical function 
$\mathfrak{P}^\circ(x,\la;q,k)$ is bounded
in terms of $x$ as
$\mathfrak{C}_+\ni \Re(x)_+\to \infty$ if and only if
$$
0< (\Re(\la)_+,\al_i^\vee)\le \Re(k_i) \for  
i=1,\ldots,n,\, \hbox{ which\ implies\ } \Re(k_\nu)> 0.
$$
If $\Re(\la)_+\in \overline{\mathfrak{C}}_+$ is allowed,
then $\mathfrak{P}^\circ_{\dagx}(x,\la;q,k)$ asymptotically
approaches a polynomial in terms of
$\{x_i\}$ of degree no greater than $n$, 
the rank of the root system.\sq
\end{corollary}

The dependence of $x$ in the right-hand side of
(\ref{growthexfinx}) is
as in the Harish-Chandra formula \cite{HC}. 
Accordingly, Corollary \ref{HELGJ}
is a $q,t$\~version of 
the description of the bounded spherical functions 
from \cite{HJ}.   

The corresponding 
degeneration of $\HH\,$ (and all related objects) is the 
procedure $q\to 1$, where we set $X_b=e^{-z_b}$ and
$z_b,\la_b$, $k$ are
considered the basic new variables upon the
degeneration. We take $-z_b$ here because the base $q$ is
smaller than $1$. The limit
of the right-hand side of (\ref{growthexfinx})  
can be readily controlled using the functional
equation for the theta-function
$\tga^\ominus$ (see below). Up to some renormalization, it 
becomes (for large $\Re(z_+)$):
$$
\hbox{Const}\,\prod_{i=1}^n\,Z_i^{(\al_i^\vee,\,\la)-k_i} 
\prod_{\al \in R_+}\, \frac{\Gamma(\la^\vee_\al)}{
\Gamma(\la^\vee_\al+k_\al)} \for Z=e^z,\, z=z_+,\la=\la_+.
$$
The factor $q^{-(\la_+,\,\rho_k)}$, which ensures
the $X\!\!\leftrightarrow\!\!\La$\~ duality of the $q,t$\~formula,
vanishes in the limit. The duality collapses under the
degeneration to the Harish-Chandra theory, however, the 
evaluation formula survives.

Technically, (\ref{growthexfinx}) matches
the growth estimates for {\em complex} Lie groups
because real Lie group result in the
terms like $\Gamma(\la_\al^\vee/2)$ in this formula, 
which is not the case.
\smallskip

The $q\to 1$ limit of the global spherical
function is convenient to describe in a somewhat
different normalization. 
Using the notations from Corollary \ref{QHARISH}, we set
\begin{align}\label{pdaggcirc}
&\mathfrak{P}^\circ_{\ddagx}(x,\la;q,k)\equal 
\frac{\mathfrak{P}^\circ(x,\la;q,k)}
{\varrho(q,t)\sigma(q^{\la_+};q,t)}\,
\frac{\tga^\ominus(q^{\la_+})\tga^\ominus(1)}
{\tga^\ominus(q^{\la_+ -\rho_k})\tga^\ominus(q^{\rho_k})}.
\end{align}
Cf. (\ref{pdagcirc}). 
Then $\mathfrak{P}^\circ_{\ddagx}(x,b-\rho_k;q,k)=
P_b(q^x;q,k)$ for all symmetric Macdonald polynomials
$P_b,b\in B_-\,$. 

Apart from the zeros of $\tga^\ominus(q^x)$,
this function is well defined for any $\la$ if all Macdonald
polynomials $\{P_b,b\in B_-\}$ exist. This condition is weaker
than the existence of all spherical polynomials
 $\{\p_b,b\in B_-\}$ we imposed above; see 
(\ref{epolexists})and (\ref{esphexists}). Moreover,
if $q^\la$ is not in the form $q^{w(b_++\rho_k)}$ for $w\in W$,
then  $\mathfrak{P}^\circ_{\ddagx}$
is well defined for {\em arbitrary} $k$ (i.e., the conditions
for $k$ necessary for the existence
of $\{P_b,b\in B_-\}$ are not needed).  

\begin{theorem}\label{PHARISH}
Provided that the Jack-Heckman-Opdam
polynomials $P\,'_b(z,\la;k)$
$=\lim_{q\to 1}P_b(e^{-z},q^\la;q,q^k)$ are well defined
for all $b\in B_-$ (a condition for $k$),
given arbitrary complex $z,\la$, the following limit
exists\,:
\begin{align}\label{qlimcirc}
&\mathfrak{P}\,'\,(z,\la;k)\ = \ \lim_{q\to 1_-}
\mathfrak{P}^\circ_{\ddagx}(e^{-z},q^\la;q,q^k).
\end{align}
This function is a $W$\~invariant solution
of the system of differential equations from \cite{HO1}
and satisfies the following conditions\,:
\begin{align*}
&\mathfrak{P}\,'_b(z,-b-\rho_k;k)\ =\ P\,'_b(z,\la;k)\for
b\in B_-\,. 
\end{align*}
Moreover, if  $\la\not\in W(B_++\rho_k)$, then the
limit $\mathfrak{P}\,'$ exists for any $k$.\sq
\end{theorem}

Here one can take complex $q=\exp(-1/a+\imath\phi)$ 
provided that $a>0, \, a\to\infty$
and $-C/a<\phi<C/a$ for a certain constant $C$. Numerical
experiments show that here $C$ can be arbitrarily large for
any given (admissible) $z,\la,k$, but we cannot justify it.

The growth estimates for $\mathfrak{P}_{\ddagx}^\circ$
read as follows:
\begin{align}\label{pdagggrow}
&\mathfrak{P}^\circ_{\ddagx}(x,\la;q,k)\ \sim\  
\frac{\tga^\ominus(q^{x_+ +\la_+-\rho_k})\tga^\ominus(1)}
{\tga^\ominus(q^{\la_+ -\rho_k})\ \tga^\ominus(q^{x_+})}\,,
\end{align}
where the asymptotic equivalence must be understood as
in Theorem \ref{GROWP} under the conditions from (i)
imposed there. The estimates are the simplest
for such normalization, since $\mathfrak{P}_{\ddagx}^\circ$
``goes through" the Macdonald polynomials.

We note that ``extending" (\ref{pdaggcirc}) and
(\ref{qlimcirc}) from the
points in the form $-b_--\rho_k$ to all $x$ is, generally,
a non-trivial problem. We involve the growth estimates.

The estimate (\ref{pdaggcirc}) becomes exactly 
$q^{-(x_+,\,\la_+-\rho_k)}$ up to a periodic function.
The latter can be readily evaluated using the following
(classical) functional equation, a progenitor
of the {\em quantum} Langlands correspondence.
The following is a variant of the formulas that
can be found in \cite{Kac}.

\begin{lemma}
Let $A$ be the lattice dual to $B$ with respect to
the standard pairing $(\,,\,)$ in $\R^n$,\, $[B:A]$\, 
the index from the theory of lattices;
for instance, $A=P^\vee $ if $B=Q$ and
$[B:A]=|P^\vee/Q|^{-1}$. Then, picking $u\in \C$ such
that $0<\Re u<\infty$,
\begin{align}\label{bfunceq}
&U^{\,x^2/2}\,\sum_{b\in B}\, X_b\, U^{\,b^2/2}\ =\ 
(\sqrt{2\pi u}\,)^{n}\,\sqrt{[B:A]}\ 
\sum_{a\in A}\, Y_a\, V^{\,a^2/2},\\
\hbox{setting\,:\ \ \ \ \,} 
&U=\exp(-\frac{1}{u}),\, V=\exp(-\frac{1}{v}) 
\for v\equal\frac{1}{4\pi^2u}\,,\,
y\equal\frac{x}{2\pi \imath u}\,,\notag
\end{align}
where $X_b=U^{x_b},\, Y_a=V^{y_a}$; for complex
$u,v$, we set $U^z=\exp(-z/u)$ and $V^z=\exp(-z/v)$.
   
Recall that $(x,x)/2=x_1^2-x_1x_2+x_2^2$ for $x_i=x_{\om_i}$
as $B=P$ in the case of $A_2$; correspondingly,
$(y,y)/2=(y_1^2+y_1y_2+y_2^2)/3$
for $A=Q$,\, $y_i=y_{\al_i}=(2x_i-x_{i\,'})/(2\pi \imath u)$,
where $\,i\,'=3-i,\, i=1,2$ for $A_2$. \sq
\end{lemma}
\smallskip

Claim (i) of Theorem \ref{GROWP}  can be naturally modified 
toward the Whittaker limiting procedure as follows.
\smallskip

\subsection{When
\texorpdfstring{{\mathversion{bold}$k\to\infty\, $}}
{t goes to 0 }}
Let us reformulate (\ref{growthek}) entirely
in terms of the function $\Psi$ from (\ref{pexla}). 
Namely, provided the conditions from Theorem \ref{GRTHM},(i), 
the limit of the function
\begin{align}\label{growthepsi}
&\Psi_{\dagx}(x,\la;q,k)\equal
\bigl(\tga^\ominus (q^{x_+ +\la_+ -\rho_k})\bigr)^{-1}\,
\Psi(X,\La;q,t)
\end{align}
exists. Similarly,
$\widetilde{\Psi}_{\dagx}(x,\la;q)\equal
\bigr(\tga^\ominus(q^{\la_+-x})\bigl)^{-1}
\,\widetilde{\Psi}(X,\La;q).$
The Whittaker limit becomes simply:
$$
\lim_{k\to \infty} 
\Psi(q^{x-\rho_k},q^\la;q,q^k)\ =\ \widetilde{\Psi}(q^x,q^\la;q).
$$
See (\ref{whitpsi}). 
\smallskip

Given real $k'_\nu \ge 0$, let us replace 
$k$ by $k+n'k'$ for $n'\in \N$ in (\ref{growthe}) and analyze the
limit 
\begin{align}\label{growthek}
&\lim_{\{\mathbf{n},n'\}\to\infty}
\Psi_{\dagx}(x+\nx,\la+\nbla;q,k+n'k').
\end{align}
In the non-simply-laced case, $n'$ can be 
treated as a $2$\~vector $\{n'_{\sht},n'_{\lng}\}$
and $n'\cdot k'$ considered instead of $n'k'$;
then both components are supposed to approach infinity
(in this paper).
\smallskip
\newtheorem{conjecture}[theorem]{Conjecture}

\begin{theorem}\label{GRTHM}
We represent $k'_\nu=u_\nu +v_\nu$ for non-negative
real $u_\nu,v_\nu$ and pick  the directions 
$\mathbf{x'},\bla$ such that
\begin{align}\label{xlarhow}
&(a)\ (\nx)_+ -n'\rho_u\,\in\, \mathfrak{C}_+\, \ni\, \Re(\la)_+
\when  \bla=0 \hbox{\ \ or\ }\\
&(b)\ (\nx)_+ -n'\rho_u\,\in\, \mathfrak{C}_+\, \ni\,  
(\nbla)_+ -n'\rho_v \when \bla\neq 0 
\notag
\end{align}
for all $\mathbf{n},n'$.
Then the limit (\ref{growthek}) exists
subject to conditions from part (i) of 
Theorem \ref{GROWP}, including the
strict positivity
requirement  $\Re(\la)_+\in \mathfrak{C}_+$.
It does not depend on
$x$ in case (a) and is a 
$x,\la$\~constant under (b). If $k'>0$ then the 
limit does not depend on $k$ too, i.e., depends only 
on $x$ for (a) and is an absolute constant for (b).
\sq
\end{theorem} 
\smallskip
The justifications are based on the formulas
for the asymptotic difference equations for the
functions under consideration. 
Theorem \ref{GROWP} corresponds to the case $k'=0$;
then the limit does depend on $k$. The rule here is
that the limit does not depend on the vectors
$x$, $\la$ or $k$ involved in the
limit, provided that the corresponding directions and the
values of the vectors which are fixed are generic.
 
The Whittaker limiting procedure can be treated as
an {\em extreme case} of the theorem as follows.
Let $k=n'k'$ assuming that $k'>0$
and  $\mathbf{n}=\{n'\}$.
We take $\la'=0$, $x'=-\rho_{k\,'}$.
Then the limit (\ref{growthek})
still exists but now it depends on $x$ (and 
depends on $\la$ too because we set $\la'=0$).
Explicitly,
\begin{align*}
\Psi_{\dagx}(x,\la;q,k)\ =\  
\bigl(\tga^\ominus(q^{\la_+ -x})\bigr)^{-1}\,
\,\Psi(q^{x+n'\,x'},q^\la; q,q^{n'\,k'}),
\end{align*}
since $(x+n'x')_+=\rho_k-x$ for sufficiently large $n'$.
Actually, we do not need $\Psi_{\dagx}$ here; the
correction factor $(\tga^\ominus(q^{\la_+ -x}))^{-1}$
does not depend on $n'$. We arrive at the procedure from
(\ref{whitgen}). 
\smallskip

We believe that the following calculation is
clarifying. Let us take {\em generic extreme} $x'$ and $\la'$ 
in (\ref{xlarhow}): 
\begin{align*}
&x'\ =\ \rho_u\,,\  \la'\ =\ \rho_v, \hbox{\ \ so\ \ } 
x'_+ +\la'_+ -\rho_k'=0.
\end{align*}
Similar to the Whittaker case, we do not need $\Psi_{\dagx}$
here.
Assuming that all $u_\nu$ and $v_\nu$ are {\em nonzero},
\begin{align}
&\lim_{n'\to\infty}\Psi(x+n'\,x',\la+n'\,\la';q,k+n'k')\notag\\
=\ 
&\sum_{b\in B_+}\, q^{\frac{(b,b)}{2}} 
\frac{X_{b}\, \La^\varsigma_{b} } {
\prod_{i=1}^n\prod_{j=1}^{(\al_i^\vee,\,b)}
(1-q_i^{j})
}\,
\label{psiuvlim}
\end{align}
for $\La^\varsigma=w_0(\La^{-1})$.
Thus, we obtain a non-constant dependence on $x$ and $\la$
here, but the output is (one of the variants of) the
multi-variable $q$\~exponential function, i.e.,
significantly simpler than the Whittaker function. 

A Whittaker variant of this calculation is
actually an extreme case of formula (\ref{growthexwla}).
It is: 
\begin{align}\label{growthexwlaext}
&\lim_{\mathbf{n}\to\infty}
\widetilde{\Psi}
(x+(\ny)_+,\la+\ny;q)\\
&=\
\sum_{b\in B_+}\, q^{\frac{(b,b)}{2}} 
\frac{\,q^{x_b-\la_b}\,} {
\prod_{i=1}^n\prod_{j=1}^{(\al_i^\vee,\,b)}
(1-q_i^{j})}\,,\notag
\end{align}
where we use the
same $\mathbf{y'}$ for $x$ and $\la$ (but in 
somewhat different way), 
assuming that $\Re(\ny)_+\in \mathfrak{C}_+$.
Note the sign of $(\ny)_+$; the growth estimates
for the $q$\~Whittaker functions considered above required
taking the direction from the negative Weyl chamber.
The proof is simple; we only need to know the leading
coefficient of $\overline{P}_b$ is $1$.
\smallskip

{\bf Discussion}.
The theorems guarantee
exponential growth (to be exact, no greater) of the function
$\mathfrak{P}^\circ$ including the boundaries of the 
domains in the theorems.

In more detail, the Gaussian-type corrections used
in the definitions of $\mathfrak{P}$\~functions
and the corresponding $\Psi$\~functions are 
not sufficient to ensure the existence of the
limits on the boundary of the domains considered
in Theorem \ref{GROWP} and \ref{GRTHM}.
Even if they are sufficient for the convergence
(as in the Whittaker case), then the limits can
depend on the initial $x,\la$. For instance, 
when $(\nx)_+,(\nbla)_+$ belong to faces of the Weyl chamber
$\overline{\mathfrak{C}}_+$, the limits
are expected to be
connected with the spherical (and Whittaker) functions
for {\em subsystems} of $R$.

The role of the condition $\Re(\la)_+\in \mathfrak{C}_+$
as $\bla=0$ is also important and not clarified
in full. As it was claimed, if
$\Re(\la)_+\not\in \mathfrak{C}_+$ 
then (\ref{growthepsi}), generally, diverges,
but the growth is polynomial.

A description of such and similar {\em
extreme} situations 
and the corresponding {\em asymptotic}
systems of  difference equations is a natural
challenge. 

Numerical experiments in the 
rank one case confirm that the convergence
condition (\ref{xlarhow}) is sharp. 
It is not clear what happens if $x'$ is taken non-proportional
to $\rho_{k\,'}$ (especially in  the non-simply-laced case
when $k=\{k_{\sht},k_{\lng}\}$). Generally,
for any $x'\in \mathfrak{C}_+$, the convergence
of $\Psi(x+n'x',\la;q,k+n'k')_{\dagx}$ 
is granted for
$0\le k'<k_o$, where $k_o=k_o(x')>0$.
What is  the formula for $k_o(x')$ and for which $x'$ the
limit exists at such extreme $k_o(x')$\,? 
\medskip

\subsection{Jackson integrals}
We are going to integrate the product of two
global spherical functions for the $\mu$\~measure
twisted by the plus-Gaussian. The previous section
guarantees that the growth of this function
in real directions is no greater than exponential. 
Due to the presence of the Gaussian,
this is sufficient to ensure the convergence of the
Jackson summations in the theorem below. 
This theorem is not from \cite{C5}, but its proof is
based on the same technique (see also \cite{C101}).

Let us fix $\xi\in \C^n$ and define the Jackson
summation as follows: 
$$\langle f\rangle_\xi\equal
|W|^{-1}\sum_{w\in W, b\in B}f(q^{w(\xi)+b}),
\where w(\xi)+b=(bw)(\!(\xi)\!).
$$ 
Here the affine action of $\hW$ from
(\ref{afaction}) is used; $f$ can be any function well defined
at the set $\{q^{w(\xi)+b}\}$. Recall that
the notation $\lan f\ran$ was used for the constant
term of a Laurent series $f$. We continue to 
assume that $|q|<1$. 

As above,  $X_{\al}(q^\xi)=q^{(\al,\xi)}, \ga(q^z)=q^{(z,z)/2}$,
$(z,z)=\sum_{i=1}^n z_iz_{\al_i}$, say,
$(z,z)/2=z_1^2-z_1z_2+z_2^2$ for $A_2$. 
For instance, 
$$\langle \ga\rangle_\xi= \sum_{a\in B} q^{(\xi+a,\xi+a)/2}=
\tga^{\ominus}(q^\xi)q^{(\xi,\xi)/2},\ 
\ga\tga^{\ominus}=\sum_{a\in B}\, a(\!(\ga)\!).
$$ 
We will constantly use that $\langle \ga\rangle_\xi$ is
periodic with respect to the substitutions
$\xi\mapsto \xi+b, b\in B$. As in \cite{C5}, let us introduce
the function $\widehat{\mu}(X;t)\equal\mu^{-1}(X;t^{-1})$,
a counterpart of $\mu$ in the theory of Jackson integration;
all $t_\al$ must be replaced by $t_\al^{-1}.$

Let us set 
$\mu_\bullet(q^{w(\xi)+b})\equal\mu(q^{w(\xi)+b})/\mu(q^\xi)$. 
Explicitly, 
using the sets $\la(bw)=\tR_+\cap (bw)^{-1}(-\tR_+)$,  
\begin{align}
& \mu_\bullet(q^{w(\xi)+b})= \!\!\!
\prod_{[\al,\nu_\al j]\in \la(bw)}
\Bigl(
\frac{
t_\al^{-1/2}-t_\al^{1/2} q^{(\al,\xi)+\nu_\al j}}{
t_\al^{1/2}-t_\al^{-1/2} q^{(\al,\xi)+\nu_\al j}
}
\Bigr)=\widehat{\mu}_\bullet(q^{w(\xi)+b})\,.
\label{muval}
\end{align}
In terms of the action 
$\hw(f)(q^z)\equal f(q^{\hw^{-1}(\!(z)\!)})$ on functions $f(q^z)$,
\begin{align}\label{lanbullet}
\lan f\mu_{\bullet}\ran_\xi\ =\ 
&\frac{\sum_{\hw\in\hW}\,\hw(f\widehat{\mu})(q^\xi)}
{|W|\,\widehat{\mu}(q^\xi)}\\
=\ 
|W|^{-1}\!\!&\sum_{w\in W, b\in B}\,f(q^{w(\xi)+b})\,
\mu_{\bullet}(q^{w(\xi)+b}).\notag
\end{align}

\begin{theorem}\label{JACKGLOB}
For arbitrary weights $\La=q^\la,\La\,'=q^{\la\,'}$,
\begin{align}\label{planjack}
(\tga^{\ominus}(q^{\rho_k}))^2\,
&\langle
\mathfrak{P}^\circ(X,\La) 
\mathfrak{P}^\circ(X^{-1},\La\,')\,\ga\mu_{\bullet}
\rangle_\xi\\ 
=\ 
\langle \ga\mu_\bullet\rangle_\xi\, 
\tga^{\ominus}_{\la}\tga^{\ominus}_{\la'}\,
\mathfrak{P}^\circ(\La,\La\,')
&\prod_{\al\in R_+}\prod_{ j=1}^{\infty}\Bigl(\frac{ 
1- q_\al^{(\rho_k,\al^\vee)+j}}{
1-t_\al^{-1}q_\al^{(\rho_k,\al^\vee)+j} }\Bigr)\,,\notag\\
\langle \ga{\mu_\bullet}\rangle_\xi=
\frac{\langle \ga\rangle_\xi}
{|W|\,\widehat{\mu}(q^\xi)}\,
&\prod_{\al\in R_+}\prod_{ j=0}^{\infty}
\Bigl(\frac{
1- t_\al^{-1}q_\al^{-(\rho_k,\al^\vee)+j}}{
1- q_\al^{-(\rho_k,\al^\vee)+j} }\Bigr).
\notag
\end{align}
\label{PLANJACK}\sq
\end{theorem}
In these formulas, $t_\nu$ are arbitrary provided the existence
of all $\{\p_{b}\}$.  The products are  
considered as the limits if 
$k_\nu\in \Z_+\setminus\{0\}$.
The normalization factor is obtained by taking
$\La=q^{-\rho_k}, \La'=q^{-\rho_k}$. Indeed,
\begin{align} 
\mathfrak{P}^\circ(X,q^{-\rho_k})\ =&\ 
\prod_{\al\in R_+}\prod_{ j=1}^{\infty}\Bigl(\frac{ 
1- q_\al^{(\rho_k,\al^\vee)+j}}{
1-t_\al^{-1}q_\al^{(\rho_k,\al^\vee)+j} }\Bigr)
\mathfrak{P}^\circ(q^{\rho_k},q^{-\rho_k}) \and\notag\\
\langle \mathfrak{P}^\circ(X,q^{-\rho_k})& 
\mathfrak{P}^\circ(X,q^{-\rho_k})
\,\ga\mu_{\bullet}
\rangle_\xi\ \notag\\  
&=\ \prod_{\al\in R_+}\prod_{ j=1}^{\infty}\Bigl(\frac{ 
1- q_\al^{(\rho_k,\al^\vee)+j}}{
1-t_\al^{-1}q_\al^{(\rho_k,\al^\vee)+j} }\Bigr)
\mathfrak{P}^\circ(q^{\rho_k},q^{-\rho_k})
\langle \ga\mu_\bullet\rangle_\xi
\label{mehjackx}
\end{align}
due to formula (\ref{hatmuxsym}).

Theorem \ref{JACKGLOB} generalizes that 
from \cite{C5} (the case of the Macdonald polynomials).
The best way of obtaining the identities from (\ref{planjack})
is via the interpretation of $\mathfrak{P}^\circ(X,\La)$ 
as the reproducing kernel of the Fourier transform, but here the 
{\em nonsymmetric setting} is more convenient. We are going to 
follow this approach in the next paper(s).

\medskip

\subsection{The special case
\texorpdfstring{{\mathversion{bold}$\xi\!=\!\!-\rho_k$}}{}}

The theory of Jackson-Gauss integrals is essentially 
algebraic, similar to that for the constant term
functional. Analytically, we need only the exponential
growth of $\mathfrak{P}$ in real directions;
(\ref{growthe}) is more than sufficient. The growth
estimates can be equally used in the theory based on 
the {\em real integration} instead of the Jackson
summation. This theory is a $q$\~generalization of
the so-called non-compact case in the harmonic analysis
on the symmetric spaces. Formulas like (\ref{planjack}) 
hold in such theory but the corresponding factors
of proportionality (generally, periodic functions in terms
of $X$ and $\La$) are not calculated so far with a
reservation about the $A_1$\~case (see \cite{C101},
Etingof's theorem).

There is a special case when (\ref{planjack}) becomes 
a straightforward {\em algebraic} exercise; it occurs
for $\xi=-\rho_k$ taken as the starting point of the 
Jackson summation. In this case,
$\mu_\bullet (q^{w(\xi)+b})$ is 
nonzero if and only if $bw=\pi_b=bu_b^{-1}$, i.e.,
at $b_\#=\pi_b(\!(-\rho_k)\!)=b-u_b^{-1}(\rho_k)$ 
in the notations from Proposition \ref{PIOM}. One has:
\begin{align}
& \mu_\bullet(q^{b_\sharp}) = 
q^{2(b_-,\rho_k)}\prod_\nu t_\nu^{l_\nu (u_b)}
\prod_{[\al,j]\in \la'(\pi_b)}
\Bigl(\frac{
1-t_\al q_\al^{(\al^\vee,\rho_k)+j}}{
1-t_\al^{-1} q_\al^{(\al^\vee,\rho_k)+j}
}
\Bigr),
\label{muvalrho}
\end{align}
where
$\la'(\pi_b)=\{\,[\al,j]\,\mid\,[-\al,\nu_\al j]
\in \la(\pi_b)\,\}$.
Then $\lan\ga\ran_{\xi}= \lan\ga\ran_{\rho_k}$ and
\begin{align}
& \langle \ga\mu_\bullet\rangle_{-\rho_k}\ =\
|W\!|^{-1}\langle \ga\rangle_{\rho_k}
\prod_{\al\in R_+}\prod_{ j=1}^{\infty}\Bigl(\frac{
1- q_\al^{(\rho_k,\al^\vee)+j}}{
1-t_\al^{-1}q_\al^{(\rho_k,\al^\vee)+j} }\Bigr).
\label{hatmu}
\end{align}
 
Formula (\ref{mehjackx})  reads as follows:
\begin{align}
&|W\!|\ (\tga^{\ominus}(q^{\rho_k}))^2\  
\langle \mathfrak{P}^\circ(X,\La) 
\mathfrak{P}^\circ(X^{-1},\La')
\,\ga\mu_{\bullet}
\rangle_{-\rho_k}\ \notag\\  
&=\ \lan \ga\ran_{\rho_k}\,
\tga^{\ominus}_\la \tga^{\ominus}_{\la'}\, 
\prod_{\al\in R_+}\prod_{ j=1}^{\infty}\Bigl(\frac{ 
1- q_\al^{(\rho_k,\al^\vee)+j}}{
1-t_\al^{-1}q_\al^{(\rho_k,\al^\vee)+j} }\Bigr)^{\!2}\
\mathfrak{P}^\circ(\La,\La'). 
\label{mehjackxx}
\end{align}

It is important to note that (\ref{mehjackxx}) is not a
{\em new} identity. To be more precise,
it formally results from the definition
of $\mathfrak{P}^\circ$, the duality of the $P$\~polynomials
and the Shintani-type relations from (\ref{hatmuxsym}). 
In a sense, the Jackson integrals trivializes at 
the special $\xi=-\rho_k$, which is analogous to the 
normalization condition in the theory of spherical functions.  
Actually relations (\ref{hatmuxsym}) were deduced in \cite{C5} 
from the general $\xi$\~theory of Jackson integration, so 
this analogy is with reservations. 
\smallskip

\subsection {Taking the limit}
Let us interpret the identity (\ref{mehjackxx}) upon
the Whittaker limit. The Jackson summation will be now
over $B$: $\lan f \ran_\diamond\equal\sum_{b\in B}f(q^{b})$;
notice that there is no $|W|$\~factor versus the
previous definition.
For instance, $\lan \ga\ran_\diamond=\tga^{\ominus}(1)$.
The corresponding 
$\overline{\mu}$\~measure is nonzero only on $B_+$:
\begin{align}\label{jackwhit}
&\overline{\mu}_\diamond(q^{b_+})\ =\ 
\prod_{i=1}^n\prod_{j=1}^{(\al_i^\vee,\,b_+)}
(1-q_i^{j})^{-1}. 
\end{align}
We come to the following ``reformulation"
of the definition of $\widetilde{\mathfrak{P}}^\circ$:
\begin{align}
&\lan \ga\ran_\diamond\,\langle
q^{(x,\la)}\, 
\widetilde{\mathfrak{P}}^\circ(X^{-1},\La')
\,\ga\overline{\mu}_{\diamond}
\rangle_{\diamond}\ \notag\\  
&=\ \tga^{\ominus}_\la \tga^{\ominus}_{\la'}\, 
\prod_{i=1}^n\prod_{ j=1}^{\infty}\Bigl(\frac{ 
1}{1-q_i^{j} }\Bigr)\
\widetilde{\mathfrak{P}}^\circ(\La,\La'). 
\label{mehjackxxx}
\end{align}
Here the Shintani-type formulas were employed.

It is instructional to obtain (\ref{mehjackxxx}) as
a Whittaker-type limit of (\ref{mehjackxx}).
We suggest the following way. 
\smallskip

First, let us
make $k$ a positive {\em integer}; to be exact, 
$k=N$ for $N=\{N_\nu\in \N\}$. Then $\rho_N\in P_+$
and, for instance, $\lan{\ga}\ran_{\la+\rho_N}= 
\lan{\ga}\ran_{\la}$, which will be used constantly.
Second, let us 
renormalize the $\mu$\~measure (\ref{muvalrho}):
\begin{align}
& \widetilde{\mu}_\bullet(q^{b_\sharp}) = 
q^{-2(b_-,\,\rho_k)}\,\mu_\bullet(q^{b_\sharp})\notag\\
=\ &\prod_\nu t_\nu^{l_\nu (u_b)}
\prod_{[\al,j]\in \la'(\pi_b)}
\Bigl(\frac{
1-t_\al q_\al^{(\al^\vee,\rho_k)+j}}{
1-t_\al^{-1} q_\al^{(\al^\vee,\rho_k)+j}
}
\Bigr).
\label{muvalrhoy}
\end{align}
The limit of $\widetilde{\mu}_\bullet(q^{b_\sharp})$ 
as $N\to\infty$ exists for any $b\in B$ and
is nonzero only for $b=b_-$. Namely,
$$
\lim_{N\to\infty} \widetilde{\mu}_\bullet(q^{b_-})=
\overline{\mu}_\diamond(q^{b_+}).
$$
Third, we will use the following property of the
{\em spherical} polynomials: 
\begin{align}\label{psphlim}
&\lim_{N\to\infty}\,\p_{b}(q^{\la-\rho_N})=q^{(\la,b_+)}
\for b=b_-,
\end{align}
which is a reformulation of (\ref{limpct}).

Forth, we observe that the condition $\rho_N\in B$
guarantees that 
$$
\lan \ga \ran_{x-\rho_{{}_N}}\ =\ 
\sum_b q^{(b+x-\rho_{{}_N},\, b+x-\rho_{{}_N})/2}
\ =\ \lan \ga \ran_x.
$$ 
for any $x$. Therefore the limiting
procedure for obtaining
$\widetilde{\mathfrak{P}}^\circ$ from 
$\mathfrak{P}^\circ$ from (\ref{whitgen})
{\em coincides} with that for 
$\widetilde{\mathfrak{P}}$ from
(\ref{whitgenmin}):
\begin{align}\label{whitgenmin1}
& \widetilde{\mathfrak{P}}^\circ(X,\La)\ =\  
\lim_{N\to \infty} q^{(x\,,\,\rho_{{}_N})}
\mathfrak{P}^\circ(q^{-\rho_{{}_N}}X,\La).
\end{align}
\smallskip

Replacing now $\la$ by $\la-\rho_N$ in
(\ref{mehjackxx}), one obtains:
\begin{align}
&\lan\ga\ran_{\diamond}^{-1}\,
(\tga^\ominus(q^{\rho_N}))^2
\langle
\mathfrak{P}^\circ(X,q^{\la-\rho_N}) 
\mathfrak{P}^\circ(X^{-1},\La')
\,\ga \mu_{\bullet}
\rangle_{-\rho_N}\ \label{mehjackxxy}\\  
=\{q^{-\frac{(\la-\rho_N)^2}{2}}&\lan\ga\ran_{\la}\} 
\,\tga^{\ominus}_{\la'}
\prod_{\al\in R_+}\prod_{ j=1}^{\infty}\Bigl(\frac{ 
1- q_\al^{(\rho_N,\,\al^\vee)+j}}{
1-t_\al^{-1}q_\al^{(\rho_N,\,\al^\vee)+j} }\Bigr)^{\!2}\,
\mathfrak{P}^\circ(q^{-(\la-\rho_N)},\La').
\notag 
\end{align}
Here $|W\!|$ is not present due to our definition
of the Jackson summation in the Whittaker case.

We can restrict ourselves only with $b=b_-$, since
the other $b$ appear in (\ref{muvalrhoy})
with strictly positive $t$\~factors 
$\prod_\nu t_\nu^{l_\nu (u_b)}$. Then the left-hand side
of (\ref{mehjackxxy}) 
{\em modulo higher powers of\, $t$\,} is as follows:\  
$LHS \hbox{\,mod\,}(t)=$
\begin{align}
&Q\,\sum_{b\in B_-}\, q^{\frac{(b_-)^2}{2}} 
\widetilde{\mu}_\bullet(b_-)
\mathfrak{P}^\circ(q^{b_-\,-\rho_{{}_N}},q^{\la-\rho_{{}_N}}) 
\{q^{(b_-,\,\rho_{{}_N})}\,
\mathfrak{P}^\circ(q^{-b_+-\rho_{{}_N}},\La')\}
\notag\\
=\
&Q\,\sum_{b\in B_-}\, q^{\frac{(b_-)^2}{2}} 
\widetilde{\mu}_\bullet(b_-)\,
\{\p(q^{\la-\rho_{{}_N}})\,\Pi\}\, 
\{q^{(-b_+,\,\rho_{{}_N})}\,
\mathfrak{P}^\circ(q^{-b_+-\rho_{{}_N}},\La')\}
\, \notag\\
&\for Q\equal \{\,q^{\rho_{{}_N}^2/2}\,
\tga^\ominus(q^{\rho_{{}_N}})/\lan \ga\ran_{\rho_{{}_N}}\,\}\,
\tga^\ominus(q^{\rho_{{}_N}})\ =\ \tga^\ominus(q^{\rho_{{}_N}})
\notag\\ 
&\and \Pi\equal 
\prod_{\al\in R_+}\prod_{ j=1}^{\infty}\frac{ 
1- q_\al^{(\rho_N,\,\al^\vee)+j}}{
1-t_\al^{-1}q_\al^{(\rho_N,\,\al^\vee)+j} }\,.\notag
\end{align}

Transforming correspondingly the right-hand side
of (\ref{mehjackxxy}), one arrives at
\begin{align} 
&RHS\ =\ q^{-\frac{\rho_{{}_N}^2}{2}}\ 
\{q^{-\frac{\la^2}{2}} \lan\ga\ran_{\la}\} 
\ \tga^{\ominus}_{\la'}\,\Pi^2\ 
\{q^{(\la,\,\rho_{{}_N})}\,
\mathfrak{P}^\circ(q^{\la-\rho_{{}_N}},\La')\}.
\notag 
\end{align}
The term  $q^{-\rho_N^2/2}$ can be moved to the LHS and
combined with $Q$, namely,
$$
q^{\rho_N^2/2}\,\tga^\ominus(q^{\rho_{{}_N}})=
\lan\ga\ran_{\rho_{{}_N}}.
$$ 
One $\Pi$ can be reduced in the LHS\,mod\,$(t)$
and the RHS. 

Then we
use (\ref{psphlim}) for
$\p(q^{\la-\rho_{{}_N}})$ and the definition of 
the Whittaker limit (\ref{whitgenmin1}) for 
$\mathfrak{P}^\circ(q^{-b_+-\rho_{{}_N}},\La')$ and
for $\mathfrak{P}^\circ(q^{\la-\rho_{{}_N}},\La')$.
Replacing (back)  
$q^{-\frac{\la^2}{2}} \lan\ga\ran_{\la}$ by
$\tga^\ominus(\la)$ and changing the summation set
in the LHS from $B_-$ to $B_+$, we eventually obtain
(\ref{mehjackxxx}).
\smallskip
 
This calculation is expected to be a sample for the
general $\xi$\~Jackson integration theory in the
Whittaker case (presumably, for the
real integration too); it will be discussed elsewhere.
We note that the term $q^{(x,\la)}$ in the
integrand of (\ref{mehjackxxx})
can be naturally combined with $\ga=q^{x^2/2}$
and ``eliminated"  upon the change
of variables $x+\la\mapsto x$. However this
substitution will change the summation
set from $B_+$ to $\la+B_+$, i.e.,
the general Jackson summation (with an arbitrary
starting vector) naturally emerges
even in the special case under consideration.
\smallskip

{\bf The extreme case.}
There is no ``natural" way to eliminate $\La,\La'$ from
(\ref{mehjackxxx}) by evaluating this formula at certain special
points. Generally, such elimination is a standard way
of discovering new identities that contain only $q$. 
Another possibility is in taking $\la\sim\la'\to \infty$  
for $\La=q^\la, \La'=q^\la$; t  Let us perform this
calculation in detail.  

We will use (\ref{growthexwla}):
\begin{align}\label{growthexwlaextu}
&\lim_{\mathbf{n}\to\infty}
\widetilde{\Psi}
(\la+(\ny)_+,\la'+\ny;q)\\
&=\
\sum_{b\in B_+}\, q^{\frac{(b,b)}{2}} 
\frac{\,q^{\la_b-\la'_b}\,} {
\prod_{i=1}^n\prod_{j=1}^{(\al_i^\vee,\,b)}
(1-q_i^{j})}\,,\notag
\end{align}
where $(\ny)_+\in \mathfrak{C}_+$.
In this limit, formula (\ref{mehjackxxx}) reads as: 
\begin{align}
&\lan \ga\ran_\diamond\,\sum_{b\in B_+}\,
q^{(b,\,\la)}\, 
\frac{\tga^\ominus(q^{\la'_+ + b})}
{\tga^\ominus(q^{\la'_+})\tga^\ominus(q^{b})}\,
\ga(q^{b})\overline{\mu}_{\diamond}(q^{b})
\notag\\  
&=\  
\sum_{b\in B_+}\, q^{\frac{(b,b)}{2}} 
\frac{\,q^{\la_b-\la'_b}\,} {
\prod_{i=1}^n\prod_{j=1}^{(\al_i^\vee,\,b)}
(1-q_i^{j})}\,,
\label{mehjackxxram}
\end{align} 
where we canceled out  
$\lan \overline{\mu}\ran=$
$\prod_{i=1}^n\prod_{ j=1}^{\infty} 
(1-q_i^{j})^{-1}$ in both sides.
Moving $b$ from the arguments of $\tga^\ominus$
and using that $\lan\ga\ran_\diamond=
\tga^\ominus(1)$, we come to an {\em identical equality}.
No new formulas appear in this way.
\smallskip

{\bf Discussion.}
We think that the growth estimates and formula
(\ref{mehjackxxx}) show great 
potential of the $q$\~theory of
Whittaker functions in harmonic analysis.
For instance, an immediate interpretation of
(\ref{mehjackxxx}) is the fact that
the global $q$\~Whittaker function multiplied
by the Gaussian is essentially {\em invariant}
with respect to the $q$\~Fourier-Jackson transform for 
the measure $\overline{\mu}_\diamond$ from
(\ref{jackwhit}), which is very much standard in
the theory of $q$\~functions. 

This paper seems a convincing 
demonstration of the key role of Shintani-type formulas
in the theory of spherical
and Whittaker functions. Interestingly, quite a 
few analytic facts (e.g, the $q$\~generalization 
of the Harish-Chandra asymptotic formula) are directly
related to these formulas. This is different
from the differential setting and makes the $q$\~theory
significantly more algebraic than the classical harmonic 
analysis on the symmetric spaces.

We would like to mention that
global spherical and Whittaker functions are expected 
to have properties similar to celebrated Ramanujan's
mock theta functions, including the theory at $|q|=1$
and certain (but not direct) counterparts
of Maas forms. To be more exact, the natural objects associated
with $q$\~spherical functions are
Maas-type  {\em theta functions}\,, which are
not holomorphic in terms of $x,\la$ but satisfy the modular
equation with respect to $q$. 

It must not be very surprising because
the {\em basic hypergeometric function}
is known to be related to (some) mock functions.
Our {\em global spherical functions} are its multi-variable
generalizations. 

\comment{
The following claim can be checked under some conditions
for $X,\la,\tilde{q}$\,; see Theorem \ref{PHARISH} for
the definition of $\mathfrak{P}^\circ_{\ddagx}$.

\begin{conjecture}\label{PHARISHROOT}
Let $\tilde{q}\,$ be unimodular. We assume that all
Macdonald polynomials $P_b$ are well defined at $\tilde q\,$ (a 
condition for $\tilde{t}= \tilde{q}^{\ k}$). Given arbitrary complex
$X,\la$, the limit exists\,:
\begin{align}\label{qlimcircroot}
&\mathfrak{P}^{\thicksim}(X,\la; k)\ = \ \lim_{q\to \,\tilde{q}_-}
\mathfrak{P}^\circ_{\ddagx}(X,q^\la;q,q^k)\,,
\end{align}
where by $q\to \,\tilde{q}_-$ we mean that $|q|<1$.
Moreover, if $X$ and $\la$ are taken ``real" with respect
to the anti-involution $f(x)\mapsto f(\bar{x}^\varsigma)$
for the complex conjugation $x\mapsto \bar{x}$ and
$\varsigma(x)=w_0(-x)$, then the limit is a real number
for $\bar{k}=k$, nonzero for generic $X,\la$.
The parameter(s) $t$ can be arbitrary here if  
$\tilde{q}^{\ \la}$ is not in the form $\tilde{q}^{\ w(b_++\rho_k)}$
for $w\in W$, $b\in B$.
\sq
\end{conjecture}
}
\medskip

\comment{
Here a reduction of the $q,t$\~Mehta-Macdonald formula 
in the Jackson setting as $t\to 0$ is used:
\begin{align}\label{barmuga}
&\lan\ga\overline{\mu}_{\diamond}\ran_{\diamond}
\ =\ \lan\ga\ran_{\diamond}
\prod_{i=1}^n\prod_{ j=1}^{\infty}(1-q_i^j)^{-1}.
\end{align} 
This relation can be deduced 
from the general theory of Jackson-type $q$\~Gauss 
integrals developed in \cite{C5,C101};  
cf. (\ref{mehtamul}). It is also a special
case of the Shintani-type formula (\ref{shintqg}):
\begin{align*}
&\sum_{b\in B_-} 
\frac{q^{(c-b,c-b)/2}\,  \overline{P}_b(\La) } {
\prod_{i=1}^n\prod_{j=1}^{(\al_i^\vee,\,b_+)}\,
(1-q_i^{j}) }=\ \tga^{\ominus}(\La)\overline{P}_c(\La)
\prod_{i=1}^n\prod_{j=1}^{\infty}\bigl(
\frac{1}{1-q_i^j }\bigr). 
\end{align*}
The demonstration goes as follows.
We employ a Whittaker-type limiting procedure 
replacing $\La$ by $q^c$ and tending  
$c\to -\infty$, say, letting $c=-2N\rho$ as $N\to \infty$.
Recall that $c\in B_-$. Then the formula becomes:
\begin{align*}
&\sum_{b\in B_-} 
\frac{q^{(b,b)/2}\, 
\{q^{-(b,c)}\,\overline{P}_b(q^c)\} } {
\prod_{i=1}^n\prod_{j=1}^{(\al_i^\vee,\,b_+)}\,
(1-q_i^{j}) }=\ \tga^{\ominus}(1)
\{q^{-(c,c)}\overline{P}_c(q^c)\}
\prod_{i=1}^n\prod_{j=1}^{\infty}\bigl(
\frac{1}{1-q_i^j }\bigr). 
\end{align*}
Here
$q^{-(b,c)} \overline{P}_b(q^c)$ can be algebraically
expanded in terms of {\em positive} powers of the 
parameter $q^{-N}$ with the constant term equals $1$
(it comes from the leading term $X_b$ of $P_b$).
This readily gives (\ref{barmuga}).
}
 

\medskip
\bibliographystyle{unsrt}

\end{document}